\documentclass[11pt]{article}                               
\usepackage{latexsym}
\usepackage{amssymb}
\input{epsf}

\def\AFOUR{%
\setlength{\textheight}{9.0in}%
\setlength{\textwidth}{5.75in}%
\setlength{\topmargin}{-0.375in}%
\hoffset=-.5in%
\renewcommand{\baselinestretch}{1.17}%
\setlength{\parskip}{6pt plus 2pt}%
}

\AFOUR                                           

\makeatletter
\def\section{\@startsection {section}{1}{\z@}{-3.5ex plus -1ex minus
 -.2ex}{2.3ex plus .2ex}{\large\bf}}
\def\subsection{\@startsection{subsection}{2}{\z@}{-3.25ex plus -1ex minus
 -.2ex}{1.5ex plus .2ex}{\normalsize\bf}}
\makeatother

\makeatletter
\@addtoreset{equation}{section}

\makeatother

\newcommand{\nc}{\newcommand}
\newcommand{\rnc}{\renewcommand}

\nc{\bea}{\begin{eqnarray}}
\nc{\eea}{\end{eqnarray}}
\nc{\be}{\bea}
\nc{\ee}{\eea}

\rnc{\a}{\alpha}
\nc{\ab}{\bar{\a}}
\nc{\ap}{\a^{+}}
\nc{\abm}{\ab^{-}}
\rnc{\b}{\beta}
\nc{\bb}{\bar{\b}}
\nc{\bbp}{\bb_{\zb}^{+}}
\nc{\bm}{\b_{z}^{-}}
\nc{\oa}{\overline{\a}}
\nc{\ob}{\overline{\b}}
\rnc{\gg}{\gamma}
\rnc{\d}{\delta}
\nc{\f}{\phi}
\nc{\fb}{\bar{\phi}}
\nc{\vf}{\varphi}
\nc{\p}{\psi}

\rnc{\c}{\chi}
\nc{\la}{\lambda}
\nc{\m}{\mu}
\nc{\n}{\nu}
\rnc{\o}{\omega}
\nc{\Om}{\Omega}
\rnc{\t}{\theta}
\nc{\eps}{\epsilon}
\rnc{\S}{\Sigma}
\nc{\F}{\Phi}

\nc{\trac}[2]{{\textstyle\frac{#1}{#2}}}

\nc{\ex}[1]{\mbox{e}^{\,\textstyle#1}}

\nc{\mat}[4]{\left(\begin{array}{cc}#1&#2\\#3&#4\end{array}\right)}

\nc{\som}[9]{\left(\begin{array}{ccc}#1&#2&#3\\#4&#5&#6\\#7&#8&#9%
\end{array}\right)}

\nc{\tr}{\mathop{\mbox{tr}}\nolimits}
\nc{\ad}{\mathop{\mbox{ad}}\nolimits}
\nc{\Tr}{\mathop{\mbox{Tr}}\nolimits}
\nc{\Det}{\mathop{\mbox{Det}}\nolimits}
\nc{\rk}{\mathop{\mbox{rk}}\nolimits}
\nc{\ra}{\rightarrow}
\nc{\Ra}{\Rightarrow}
\nc{\LRa}{\Leftrightarrow}
\nc{\ot}{\otimes}
\rnc{\ss}{\subset}
\nc{\nul}{\noindent\underline}
\nc{\non}{\nonumber\\}

\nc{\subs}[1]{{\vspace*{0.5cm}}%
{\noindent\underline{#1}}{\addcontentsline{toc}{subsection}{#1}}%
{\vspace*{0.3cm}}}

\nc{\zb}{\bar{z}}
\rnc{\lg}{\frak{g}}
\nc{\lt}{\frak{t}}
\nc{\lk}{\frak{k}}
\nc{\lh}{\frak{h}}
\nc{\pik}{\Pi_{\lk}}
\nc{\pip}{\Pi_{+}}
\nc{\pim}{\Pi_{-}}
\nc{\pih}{\Pi_{\lh}}
\nc{\jz}{J_{z}}
\nc{\jzh}{\jz^{\lh}}
\nc{\jzp}{\jz^{+}}
\nc{\jzm}{\jz^{-}}
\nc{\del}{\partial}
\nc{\dz}{\del_{z}}
\nc{\dzb}{\del_{\bar{z}}}
\nc{\az}{A_{z}}
\nc{\azb}{A_{\bar{z}}}
\nc{\g}{g^{-1}}
\nc{\dw}{\Delta_{W}}
\nc{\Ad}{{\mbox{Ad}}}
\nc{\ks}{Ka\-za\-ma-\-Su\-zu\-ki}
\nc{\KS}{\ks}
\nc{\ksm}{\ks\ model}
\rnc{\AA}{{\Bbb A}}
\nc{\BB}{{\Bbb B}}
\nc{\CC}{{\Bbb C}}
\nc{\PP}{{\Bbb P}}
\nc{\cpm}{\CC\PP(m)}
\nc{\cpn}{\CC\PP(n)}
\nc{\cp}[1]{\CC\PP(#1)}
\nc{\gmn}{G(m,m+n)}
\nc{\gmnk}{\gmn_{k}}
\nc{\cO}{{\cal O}}
\nc{\bcO}{\bar{\cO}}
\nc{\bO}{\bar{O}}
\nc{\oQ}{\overline{Q}}

\nc{\qed}{\indent \hspace{14cm} $\Box$}

\newtheorem{theorem}{Theorem}[section]
\newtheorem{corollary}{Corollary}[theorem]
\newtheorem{proposition}[theorem]{Proposition}
\newtheorem{definition}{Definition}[section]
\newtheorem{remark}{Remark}[section]
\newtheorem{lemma}{Lemma}[section]

\begin{document}
\global\parskip=4pt
\makeatletter
\begin{center}
\vskip .5in
{\LARGE\bf 
Skew Invariant Theory of Symplectic Groups, Pluri-Hodge Groups
and 3-Manifold Invariants}
\vskip .5in
{\bf George Thompson}\\
{\bf Abdus Salam I.C.T.P.}\\
{\bf Strada Costiera 11, Trieste, Italy}\\
{\bf thompson@ictp.it}
\end{center}
\makeatother


\begin{small}
\end{small}

\setcounter{footnote}{0}

\section*{Introduction}

This paper is concerned with four subjects. The second and third of
these appear to be, at first sight, quite disparate. The fourth,
however, shows how they are, somewhat surprisingly, related.

The first has to do with skew invariant theory. Let $V$ be a finite
dimensional complex vector space with a skew symmetric non-degenerate form and
let $\mathrm{Sp}(V)$ be the corresponding symplectic group. For a
complex vector space $W$ we consider the action of $\mathrm{Sp}(V)$ on the
exterior algebra $\Lambda(V \otimes W)$ (the action being trivial on
$W$). In Theorem \ref{fund} we give generators and relations for the
algebra of $\mathrm{Sp}(V)$ invariants for this action (``Fundamental
theorem for skew invariant theory''). The relations are the so-called
``$P_{n}$'' relations which appear in the study of certain 3-manifold
invariants.

Next we show that there is a natural action of $\mathrm{Sp}(g)$, where
$g= \dim_{{\mathbb{C}}} W$, on $\Lambda(V \otimes W)$, by considering
  the spin representation of the orthogonal algebra associated to $V
  \otimes (W \oplus W^{*})$. The action of $\mathrm{Sp}(V)\times
  \mathrm{Sp}(g)$ on $\Lambda(V \otimes W)$ is multiplicity free and
  we determine the highest
  weights of the representations which occur, Theorem \ref{decomp}.

The second subject is concerned with the study of Dolbeault cohomology
groups with values in specific vector bundles. Let $X$ be a compact
closed complex
manifold. Let $\Omega^{p}_{X}$ be the sheaf of holomorphic $p$-forms
on $X$. The cohomology groups of interest are
$$
\mathrm{H}^{q}(X,\Omega_{X}^{p_{1}}\otimes \dots
\otimes \Omega_{X}^{p_{g}}).
$$

These cohomology groups are not exotic. For example,  the
pluri-canonical sections, $\mathrm{H}^{0}(X, K_{X}^{\otimes m})$, are
well studied. In keeping with this terminology I call these
cohomology groups, pluri-Hodge groups.

Some of the motivation for studying pluri-Hodge groups, when $X$
is holomorphic symplectic, comes from physics. In section
\ref{chernnumbers} I
give an alternative motivation. Namely one can show that knowing the
dimensions of these
groups (actually one needs somewhat less) one knows all the Chern
numbers of $X$. Theorem \ref{cstring} gives a precise formula for the
Chern numbers in terms of the dimensions of the pluri-Hodge
groups. Along the way a pluri-$\chi_{y}$ genus is defined.

In section \ref{holsympman} one specializes to holomorphic symplectic
manifolds. In this setting physics suggests that there is
an $Sp(g)$ action on the direct sum of these pluri Hodge spaces. This
is a straightforward 
consequence of the preparational material in section \ref{sreps}. The
case $g=1$ has been studied by Fujiki \cite{F}. The
multiplicities of the irreducible components of the total representations are
determined in Theorem \ref{multi}. Also in this
section Theorem \ref{echi} shows that the graded trace of an $Sp(g)$ element is
essentially (up to normalization) the pluri $\chi_{y}$ genus. This
generalizes the situation for the usual $\chi_{y}$ genus.

For a hyper-K\"{a}hler manifold $X$ there are good reasons (see below) 
to believe that the dimensions of the pluri-Hodge groups do not depend on the
intrinsic $\mathbb{P}^{1}$ of complex structures associated with $X$. In Theorem
\ref{k3constant} the dimensions of the pluri-Hodge
groups are determined in the case of K3 surfaces and this shows that
indeed the dimensions do not depend on the specific K3 surface.

Obviously, one unifying theme between the two subjects is the group
$Sp(g)$. Surprisingly, though, another thread through these subjects
has to do with 
a class of 3-manifold invariants. A special case of Theorem
\ref{decomp} tells us that the space of $\mathrm{Sp}(V)$ invariants
carries an irreducible representation of $\mathrm{Sp}(g)$ (of type
$(n, \dots , n)$ where $2n = \dim_{\mathbb{C}}V$) while the
fundamental theorem of skew invariants Theorem \ref{fund} says that
the space of invariants is determined by the $P_{n+1}$ relations. 
On the other hand the $P_{n}$ relations appear in the study of certain
3-manifold invariants. 

For 3-manifolds without boundary these are the LMO \cite{LMO}
invariants and for 3-manifolds with boundary they are the MO \cite{MO}
invariants. In these theories the $P_{n}$ relations have been
implemented on spaces of chord diagrams
$\mathcal{A}(\Gamma_{g})$. There is an equivalent formulation, for the
MO invariants,
employing spaces of uni-trivalent graphs $\mathcal{B}_{g}$ introduced
by J. Sawon \cite{justinhk}. His formulation is used, in Theorem
\ref{Brep} to show that the
resulting quotient spaces carry representations of symplectic
groups. Theorem \ref{chordspg} determines the irreducible
representations that occur in the simplest case.

We come to the last topic. There is a path integral formulation of LMO
type invariants due 
to Rozansky and Witten \cite{RW}. For closed 3-manifolds there is
strong evidence that there is a weight system which when applied to
the LMO invariants produces the Rozansky-Witten invariants
\cite{HT}. The 
pluri-Hodge groups for a hyper-K\"{a}hler manifold make an
appearance in the Rozansky-Witten theory. This theory assigns to a
connected genus $g$ surface $\Sigma_{g}$ the Hilbert space of states,
$$
\mathcal{H}_{g}(X) = \bigoplus_{q, p_{1}, \dots , p_{g}} \mathrm{H}^{q}(X,
\Omega^{p_{1}}_{X} \otimes \dots \otimes \Omega^{p_{g}}_{X} ).
$$
Sawon showed that there is a weight system which is a homomorphism of
vector spaces $W_{X}:\mathcal{B}_{g} \rightarrow
\mathcal{H}_{g}(X)$ (this weight system appears also in \cite{Tunpub} in a more
limited situation). Theorem \ref{commspg} states that $W_{X}$
preserves the symplectic group actions on the spaces.  In this way one
explains the ubiquitous appearance of the representation theory of the
symplectic group. The Rozansky-Witten
theory then acts as a unifying thread for the symplectic group actions
in this paper. 

To finish this introduction it should be pointed
out how the RW theory gives rise to some expectations that are
expressed in the text. Proposition \ref{half}
states that the homomorphism $W_{X}$ has as image a much smaller space than
$\mathcal{H}_{g}(X)$. However, the generalized $\chi_{y}(X)$ genera for
holomorphic symplectic manifolds also arise in the RW theory for
certain 3-manifolds that are mapping tori \cite{T}. These `see' a lot
more of the pluri-Hodge groups than allowed by the proposition. This
should help in formulating a TQFT for these theories. For
closed 3-manifolds 
the RW theory can be expressed 
in terms of a given hyper-K\"{a}hler manifold $X$ without reference to
a preferred complex structure in the $\mathbb{P}^{1}$ of complex
structures. Since in the physical theory the RW invariant arises as
pairings of vectors in $\mathcal{H}_{g}(X)$ these pairings should not
depend on the choice of complex structure implicit in $\mathcal{H}_{g}(X)$. With
this in mind it is natural to expect that the dimensions of the
pluri-Hodge groups are deformation invariants.

{\bf Acknowledgements:} I have discussed some of these things with
Matthias Blau, I am grateful to him for listening and for reminding me
of some representation theory facts. The proof of
Theorem \ref{fund} presented here is due to C. De Concini and
C. Processi, I would like to express my thanks to them for providing
this as my earlier proof was less intrinsic and more
computational. Thanks go also to the referee for suggesting that the
Rozansky-Witten story should be included and for a myriad of
improvements to the text. All the ideas here have been 
developed over a long period with M.S. Narasimhan. He kindly
gave much of his time, a great
deal of assistance and insisted that I collect the pieces together. It
is a special pleasure to thank him for all of these things and more.

\section{Special Representations of Products of Symplectic Groups}\label{sreps}

Let $V$ be a symplectic vector space of complex dimension $2n$ with
symplectic form $\eps$. We
will define in this section a natural representation of $Sp(g)$ (with
$\mathrm{Sp}(1) = \mathrm{SL}(2, \mathbb{C})$) on the
exterior algebra of $V \otimes \mathbb{C}^{g}$ and decompose this
space under the action of $Sp(V) \times Sp(g)$. For convenience of
notation let $W$ be a vector space of complex dimension $g$. We have
on $W \oplus W^{*}$ a natural symplectic form, $\alpha$, defined by
$$
\alpha\left( (w_{1},w_{1}^{*}), (w_{2},w_{2}^{*})\right) =
w_{2}^{*}(w_{1}) - w_{1}^{*}(w_{2}).
$$
Now $h=\eps \otimes \alpha$ defines a non-degenerate symmetric bi-linear
form on
$H=V \otimes (W \oplus W^{*})$ and let $o(H,h)$ be the Lie algebra 
of the corresponding orthogonal group. Note that $Sp(V) \times Sp(W)$
leaves $h$
invariant. It is now not too surprising that the spin representation
makes an appearance.

Observe that $V \otimes W$ is a maximal isotropic space of $H$ with respect
to $h$ since $\alpha$ on $W$ vanishes. Thus we have a spin
representation of $o(H,h)$ on $\Lambda (V \otimes W)$ \cite{B}. By restriction
we get a representation of the product of the Lie algebras of the
symplectic groups and, as these groups are simply connected, we indeed
get a representation of the product group. Furthermore, we observe for
later use, that the Lie algebras of $Sp(V)$ and
$Sp(W \oplus W^{*})$ form a dual reductive pair in the Lie algebra of
$O(H)$ in the sense of Howe \cite{howe}.

\subsection{Spin Representations}

The following comes from \cite{B} (see in particular Lemma 1 p. 195).
Let $(H,h)$ be an even dimensional orthogonal vector space which can
be decomposed as 
$$H= F \oplus F^{*}$$
for $F$ a finite dimensional vector space over $\mathbb{C}$ and
$F^{*}$ its dual. The symmetric bilinear form $h$ is given by
$h(x_{1}+x_{1}^{*}, y_{1}+ y_{1}^{*}) = \langle y_{1},
x_{1}^{*}\rangle + \langle x_{1}, y_{1}^{*}\rangle$ where $x_{1}, \;
y_{1} \in F$, $x_{1}^{*}, \; y_{1}^{*} \in F^{*}$ and $\langle \, , \,
\rangle$ denotes the pairing between $F$ and $F^{*}$. Given $e_{a}$ as
a basis for $F$ and $\tilde{e}^{b}$ a dual basis for $F^{*}$, the
symmetric form is such that
$$
h(e_{a},e_{b})=0, \;\;\; h(e_{a},\tilde{e}^{b}) = \delta_{a}^{b},
\;\;\;  h(\tilde{e}^{a},\tilde{e}^{b})=0.
$$ 

Denote by $C(Q)$ the Clifford algebra of $H$ with respect to the
quadratic form $Q$ defined by $Q(v) = h(v,v)/2$. We have a canonical
map $f_{0}: o(H,h) \rightarrow C(Q)$ defined as the composition of
maps
$$
o(H,h) \rightarrow \mathrm{End}(H) \rightarrow H \otimes H^{*}
\rightarrow H \otimes H \rightarrow C(Q).
$$
Here the last map is given by multiplication in the Clifford algebra,
and the third map by using the canonical isomorphism of $H^{*}$ with
$H$ given by the bilinear form $h$. Put $f= f_{0}/2$, then $f$ is a
Lie algebra homomorphism of $o(H,h)$ into the Lie algebra
corresponding to the associative algebra $C(Q)$. Indeed, for $X\in o(H,h)$,
\be
f(X) = \frac{1}{2} \sum_{a} \left[ (Xe_{a})\otimes \tilde{e}^{a} +
  (X\tilde{e}^{a} ) \otimes e_{a} \right].\label{f}
\ee

Let $\Lambda F $ be the exterior
algebra of $F$ and $
\Lambda F^{*} $ be the exterior algebra of $F^{*}$. Let $\lambda: C(Q)
\rightarrow \mathrm{End}(\Lambda F)$ and $\lambda^{*}: C(Q) \rightarrow
\mathrm{End}(\Lambda F^{*})$ be the two spin representations of
$C(Q)$. Concretely, on $\Lambda F$ (resp. $\Lambda F^{*}$) the
defining property of $\lambda$ (resp. $\lambda^{*}$) is to send
$\tilde{e}^{b}$ (resp. $e_{b}$) to interior multiplication,
$\iota_{e_{b}}$ (resp. $\iota_{\tilde{e}^{b}}$), while
$e_{a}$ (resp. $\tilde{e}^{a}$) is understood to be exterior
multiplication (wedging on the
left), these operations then give spin representations since, on the
respective spaces, one has
$$
\{ i_{e_{b}}, e_{a} \} = \delta_{ab}, \;\;\; \{ i_{\tilde{e}^{b}},
\tilde{e}^{a} \} = \delta_{ab}.
$$

\subsection{Restriction of Spin to $\mathrm{Sp}(g)$}

The
spin representations of $o(H,h)$ are $\rho:= \lambda .f:
o(H,h)\rightarrow \mathrm{End}(\Lambda F)$ and $\rho^{*}:= \lambda^{*} .f:
o(H,h)\rightarrow \mathrm{End}(\Lambda F^{*})$. Note that the Lie
algebra of the orthogonal group is
$$
o(H,h)= S^{2}(F\oplus F^{*})= S^{2}F\oplus S^{2}F^{*} \oplus
F\otimes F^{*}.
$$

Let $S: F \rightarrow
F$ be a linear map and let $\widetilde{S}: F\oplus F^{*}
\rightarrow F\oplus F^{*}$ be defined by $\widetilde{S} = (S, - S^{T})$
then $\widetilde{S} \in o(H,h)$. In particular $\widetilde{S} \in
F\otimes F^{*}$.

Furthermore,
$$
f(\widetilde{S}) = \frac{1}{2}\sum_{a}\left[ (Se_{a})\otimes
  \tilde{e}^{a} - (S^{T}\tilde{e}^{a})\otimes e_{a} \right].
$$
We have
\begin{lemma}\label{lemma1}{\rm{Let $S: F \rightarrow F$ be a linear map and
      $\widetilde{S} = (S, - S^{T})$. Then
\begin{enumerate}
\item $\rho(\widetilde{S}) =  (S) - \frac{1}{2}\Tr{(S)}. \mathbb{I}$
\item $\rho^{*}(\widetilde{S}) =  (-S^{T}) +
\frac{1}{2}\Tr{(S)}. \mathbb{I} $
\end{enumerate}
where (S) denotes the extension to the exterior algebra $\Lambda F$ as
derivation and $(-S^{T})$ denotes the unique derivation of $\Lambda F^{*}$
which coincides with $(-S^{T})$ on $F^{*}$.}}
\end{lemma}
{\bf Proof:} We show the first statement, the proof of the second is
similar. In $C(Q)$ we have $\tilde{e}^{a} \otimes e_{a} =
\mathbb{I} - e_{a}\otimes \tilde{e}^{a}$, and after applying $\lambda$
we understand that $\tilde{e}^{a} = \iota_{e_{a}}$ etc. so that,
$$
\rho (\widetilde{S}) = \sum_{a}\left[(Se_{a}).\iota_{e_{a}} -
  \frac{1}{2} \Tr{(S)}.\mathbb{I}\right] 
$$
\qed

Note that $Gl(W^{*})$ has a natural action, $\rho_{0}$, on $\Lambda (V
\otimes W)$ with its
action on $V \otimes W$ being $\mathrm{Id_{V}} \otimes (M^{T})^{-1}$
with $ M \in
Gl(W^{*})$. On the other hand $Gl(W^{*})$ is a subgroup of $Sp(W
\oplus W^{*}) \equiv Sp(g)$. However, the restriction to $Gl(W^{*})$
of the action of 
$Sp(g)$ on $\Lambda ( V \otimes W)$ is the natural action twisted by a
character of
$Gl(W^{*})$. More precisely we have, on applying Lemma \ref{lemma1}, 
\begin{proposition}\label{central}{\rm{Let $\rho$ be the restriction
      of the spin
      representation to $Gl(W^{*})$ and $\rho_{0}$ be the natural
      representation on $\Lambda ( V \otimes W)$, $M \in
Gl(W^{*})$ then we have
$$
\rho(M) = (\det{M})^{n}\, \rho_{0}(M).
$$}}
\end{proposition}
{\bf Proof:} Apply Lemma \ref{lemma1} with $F = V \otimes W$, and $S=
\mathrm{Id_{V}} \otimes -S_{0}^{T}$ to obtain the formula of
interest at the the Lie algebra level with $S_{0} \in
gl(W^{*})$. Exponentiation yields the proposition.

\qed

\subsection{Explicit Formulae}\label{explicit}

By considering $X \in S^{2}F$ and $X\in S^{2}F^{*}$ as well one can be more
explicit about the representations. Any element of ${\rm sp}(g)$, in
the defining representation, can be written as a block matrix
\be
\left( \begin{array}{cc}
\mathbf{a} & \mathbf{b} \\
\mathbf{c} & -\mathbf{a}^{T}
\end{array} \right) \label{sp1}
\ee
of $g\times g$ matrices, where $\mathbf{a} \in {\rm gl}(g)$ and
$\mathbf{b}$ and $\mathbf{c}$ are
symmetric matrices. Let $E_{ij}$ be the $g \times g$ matrix whose only
non-zero entry is the one at the i-th row and j-th column and that
entry is 1. The identification to be made is the following
$$
h^{i}_{\;\; j} = \left( \begin{array}{cc}
E_{ji} & 0 \\
0 & -E_{ij} 
\end{array} \right) -n \delta^{i}_{\; j},
$$
with the shift by the identity matrix taking into account the
twisting by a character, and 
$$ L_{ij} = \left( \begin{array}{cc}
 0& \frac{1}{2}(E_{ij}+ E_{ji}) \\
0 & 0
\end{array} \right) ,\; \; \Lambda^{ij} = \left( \begin{array}{cc}
 0& 0 \\
\frac{1}{2}(E_{ij}+ E_{ji}) & 0
\end{array} \right).
$$
These clearly span the matrices of the type (\ref{sp1})

Set $V
\otimes \mathbb{C}^{g} = V_{1} \oplus \dots
\oplus V_{g}$ where all of the spaces $V_{i} = V$. Introduce on each
$V_{i}$ a basis of $\Lambda^{1} V_{i}$, $v^{I}_{i}$, where $I =1 ,
\dots , \; 2n$ and $i =1, \dots , \; g$. Introduce the interior
multiplication, $i^{j}_{J}$, with respect to the previous basis
$$
i^{j}_{J} \; v^{I}_{i} = \delta^{j}_{\; i} \, \delta^{I}_{\; J}
$$
and let the symplectic form, $\eps$, on $V$ be given by $\eps=
\frac{1}{2} \sum_{\{I, \, J\}}\eps_{IJ}\, v^{I} \wedge v^{J}$. Since the
symplectic form is non-degenerate there exists an inverse matrix,
$\eps^{IJ}$, such that
$$
\sum_{J=1}^{2n}\eps^{IJ}\eps_{JK} = \delta^{I}_{\; K}.
$$
Define the following operations on $\Lambda(V \otimes W)$, let
$\omega \in \Lambda(V \otimes W)$,
\begin{definition}{\rm{
$$
L_{ij}(\omega)  =   \frac{1}{2} \, \sum_{I,J}\eps_{IJ}v^{I}_{i}\wedge
v^{J}_{j}\, \wedge \omega , 
\;\;\; \Lambda^{ij}(\omega)  =  \frac{1}{2} \, \sum_{I,J} \eps^{IJ} \,
i^{i}_{I} \; 
i^{j}_{J}( \omega), 
\;\;\; h^{i}_{\; j} (\omega)  = \sum_{J} v^{J}_{j}\wedge i^{i}_{J} (\omega) 
$$
}}
\end{definition}
Note that this definition is actually an application of (\ref{f}) on the
appropriate Lie algebra elements with $\mathrm{sp}(g) \hookrightarrow
o(H,h)$.

\begin{proposition}\label{spalg1}{\rm{The operators, $L_{ij}$,
      $\Lambda^{ij}$ and 
      $h^{i}_{j}$ satisfy the relations
\bea
\left[ \Lambda^{ij}, L_{kl} \right] & = & \frac{n}{2}\left(
\delta_{l}^{i} \delta_{k}^{j} + \delta_{k}^{i} \delta_{l}^{j} \right)
- \frac{1}{4}
\left(\delta_{k}^{i} h_{l}^{j} + \delta_{l}^{i} h_{k}^{j} +
\delta_{k}^{j} h^{i}_{l} + \delta_{l}^{j} h^{i}_{k} \right) \nonumber \\
\left[ h^{k}_{\; \; l} \, , L_{ij} \right] & = & \delta_{i}^{k}L_{lj} +
\delta_{j}^{k} L_{il} , \;\;\; 
\left[ h^{k}_{\; \; l} \, , \Lambda^{ij} \right] =  - \delta^{i}_{l}
\Lambda^{kj} - \delta^{j}_{l}\Lambda^{ik}  \nonumber \\
\left[ h^{i}_{\;\; j}\, , h^{k}_{\;\; l} \right] & = & \delta^{i}_{l}
h^{k}_{\;\; j} - \delta^{k}_{j} h^{i}_{\;\; l}, \;\;\; 
\left[ \Lambda^{ij}, \Lambda^{kl} \right]  =  0, \;\;\; 
\left[ L_{ij}, L_{kl}^{ } \right]  =  0 \nonumber
\eea
and this is a representation of the Lie algebra ${\rm sp}(g)$.}}
\end{proposition}

{\bf Proof:} That the generators satisfy the algebra and that this is
$\mathrm{sp}(g)$ follows from what
has been said. Alternatively the algebra can be checked by a short
calculation and 
can be put in standard form by sending $h^{i}_{\; j} \rightarrow
h^{i}_{\; j} + n \delta^{i}_{\; j}$. \qed

Let $\alpha \in \Lambda^{p_{1}}V \otimes \dots \otimes \Lambda^{p_{g}}V $ and
make the dependence on the basis of $\Lambda^{p_{1}}V \otimes \dots
\otimes \Lambda^{p_{g}}V $ explicit as $\alpha(v_{i}^{I})$. Also set
$|p| = \sum_{i=1}^{g} p_{i}$.

The following theorem allows us to write the $Sp(g)$
action explicitly, and represents the exponentiation alluded to in Proposition
\ref{central}.
\begin{theorem}\label{rep1}{\rm{The following vector spaces are
      representations of ${\rm Sp}(g)$ 
\bea
{\mathcal H}_{+}(V) &=& \bigoplus_{\sum p_{i}= 0\; {\rm mod}\; 2}
\Lambda^{p_{1}}V \otimes \dots \otimes \Lambda^{p_{g}}V , \;\;\;
\mathrm{and}\nonumber \\ 
{\mathcal H}_{-}(V) & =& \bigoplus_{\sum p_{i}= 1\; {\rm mod}\; 2}
\Lambda^{p_{1}}V \otimes \dots \otimes \Lambda^{p_{g}}V .\nonumber
\eea
In particular for $U \in {\rm Sp}(g)$ and $\alpha \in
\Lambda^{p_{1}}V \otimes \dots \otimes \Lambda^{p_{g}}V$ the action is given by
$$
\alpha(v_{i}^{I}).U = (\Det{\mathbf{A}})^{n}\alpha(\eps^{IJ}\mathbf{C}.i_{J} +
\mathbf{D}.v^{I})\, \exp{\left(-\Tr (L
  . (\mathbf{A}^{-1}.\mathbf{B}))\right) },
$$
where
$$
U = \left( \begin{array}{cc}
\mathbf{A} & \mathbf B\\
\mathbf{C} & \mathbf{D}
\end{array}\right).
$$ }}
\end{theorem}
{\bf Proof:} By the previous proposition we know that the space
$$
\mathcal{H}(V) = {\mathcal H}_{+}(V) + {\mathcal H}_{-}(V)
$$
furnishes a representation space for the Lie algebra $sp(g)$. Indeed,
as the generators $L$, $\Lambda$ and $h$ all change $|p|$
by $0 \; \mathrm{mod} \; 2$ we see that the individual spaces, ${\mathcal
  H}_{+}(V)$ and ${\mathcal H}_{-}(V)$, form
representation spaces for the Lie algebra. Since the symplectic group is
simply connected we know that one can obtain a representation from the
algebra by exponentiation. However, this is not to say that the action
of the group is easy to write down explicitly. That the formula
advocated is a representation of $Sp(g)$ can be shown by
straightforward but rather tedious algebra which the
reader is spared. However, it is possible to easily establish its validity on
the subgroup of $Sp(g)$ where the bottom left hand block ($\mathbf{C}$)
is zero. Let $U_{1}$ and $U_{2}$ be elements of the subgroup, then
$$
(\alpha.U_{1}).U_{2} =
(\Det{\mathbf{A}_{1}.\mathbf{A}_{2}})^{n}\alpha(
 \mathbf{D}_{1}.\mathbf{D}_{2}.v^{I})\, \exp{\left(-\Tr (L
   . [\mathbf{D}_{2}^{T}. \mathbf{A}^{-1}_{1}.\mathbf{B}_{1}.
     \mathbf{D}_{2} +   \mathbf{A}^{-1}_{2}.\mathbf{B}_{2}])  \right)  }
$$
The product $U_{3}= U_{1}U_{2}$ is,
$$
U_{3} = \left( \begin{array}{cc}
\mathbf{A}_{3} & \mathbf{B}_{3} \\
\mathbf{0} & \mathbf{D}_{3}
\end{array} \right) =
\left( \begin{array}{cc}
\mathbf{A}_{1}.\mathbf{A}_{2} & \mathbf{A}_{1}.\mathbf{B}_{2} +
\mathbf{B}_{1}. \mathbf{D}_{2}\\
\mathbf{0} & \mathbf{D}_{1}.\mathbf{D}_{2}
\end{array} \right),
$$
and we note that $\mathbf{A}_{3}^{-1}.\mathbf{B}_{3} =
\mathbf{A}_{2}^{-1}.\mathbf{B}_{2} +
\mathbf{A}_{2}^{-1}.\mathbf{A}_{1}^{-1}.
\mathbf{B}_{1}. \mathbf{D}_{2}$. However, for $U_{2}$ to be an element
of $Sp(g)$, we must have that $\mathbf{D}_{2}^{T}=\mathbf{A}_{2}^{-1}$
so that we have shown
$(\alpha.U_{1}).U_{2} = \alpha.(U_{1}.U_{2})$ as required. 

The action of the group on $\alpha$ also preserves
$(-1)^{|p|}$ since all the terms appearing in the formula
shift $|p|$ by multiples of 2 and so preserves the
splitting of $\mathcal{H}(V)$.\\
\qed

\begin{remark}{\rm{The action of minus the identity,
      $-\mathbb{I}$ the generator of the centre of $Sp(g)$, is as
      multiplication by $(-1)^{ng-|p|}$. Hence, it is
      $Sp(g)/\pm \mathbb{I}$ which acts on $\mathcal{H}_{+}(V)$ when
      $ng=0 \; \mathrm{mod}\, 2$ and on $\mathcal{H}_{-}(V)$
      when $ng=1 \; \mathrm{mod}\, 2$.}}
\end{remark}

\section{Skew Invariant Theory}

Let $V$ and $W$ be finite dimensional vector spaces of dimension $m$
and $g$ over $\mathbb{C}$ respectively.
\begin{proposition}\label{schur}{\rm{Under the action of $GL(V)\times
      GL(W)$ the 
      exterior algebra $\Lambda ( V \otimes W)$ decomposes as
$$
\Lambda ( V \otimes W) = \oplus_{\lambda} \left(S_{\lambda}(V) \otimes
S_{\tilde{\lambda}}(W)\right)
$$
where $S_{\lambda}$ (respectively $S_{\tilde{\lambda}}$) is the Schur
  functor corresponding to the Young diagram $\lambda$ (respectively
  to the Schur functor of the dual partition $\tilde{\lambda}$) and
    the partition $\lambda$ has at most $m$ columns and
    $\tilde{\lambda}$ has at
    most $g$ rows.}}
\end{proposition}
{\bf Proof:} See for example Theorem (8.4.1) in \cite{procesi}.

\qed

\begin{proposition}\label{invs}{\rm{Suppose that $V$ is of complex
      dimension $2n$
      and comes equipped with a non-degenerate skew symmetric form.
      Denote by $Sp(V)$ the corresponding symplectic group, then, the
      space $\mathcal{I}$ of $Sp(V)$ invariants in $\Lambda ( V
      \otimes W)$ has the following decomposition under $GL(W)$. 
$$
\mathcal{I}= \oplus S_{\mu}(W)
$$
where
\begin{enumerate}
\item each row of $\mu$ has an even number of elements
\item the
first row has $\leq 2n$ elements and 
\item the number of rows is at most $g$.
\end{enumerate}}}
\end{proposition}
{\bf Proof:} For a partition $\lambda$, the space of $Sp(V)$
invariants in $S_{\lambda}(V)$, is denoted $S_{\lambda}(V)^{Sp(V)}$. We have,
using proposition \ref{schur},
$$
\mathcal{I}=\left( \oplus S_{\lambda}(V) \otimes
S_{\tilde{\lambda}}(W)\right)^{Sp(V)} =  \oplus\left(
S_{\lambda}(V)^{Sp(V)}  \otimes
S_{\tilde{\lambda}}(W)\right).
$$
On the other hand, (see the following remark)
\be
\dim{S_{\lambda}(V)^{Sp(V)}} = \left\{ \begin{array}{l}
1, \; \mathrm{if \; each \; column\; of\; \lambda \; has\; an\; even \;
  number\; of\; elements}\\
0,\; \mathrm{otherwise}
\end{array} \right.\label{Decomp}
\ee
Now this proposition follows from the previous one.

\qed

\begin{remark}{\rm{While (\ref{Decomp}) is apparent I have not been
      able to find a reference for it. The result follows from the
      restriction rules of \cite{kt}, page 443, together with the
      conditions on the Littlewood-Richardson coefficients given in
      \cite{m} page 142.}}
\end{remark}

\begin{theorem}\label{fund}{\rm{({\bf The fundamental theorem of skew
        invariant theory})

The (commutative) algebra $\mathcal{I}$ of $Sp(V)$ invariants in
$\Lambda ( V \otimes W)$ has the following presentation. Let
$S(S^{2}W)$ denote the symmetric algebra of the second symmetric power
  of $W$. Consider the natural embedding, $P_{n+1}$, of $S^{2n+2}(W)$
  into $S(S^{2}W)$ given by the inclusion $S^{2n+2}(W) \hookrightarrow
  S^{n+1}(S^{2}W)$. See for example (4.3.3) in \cite{GW}. Then
  $\mathcal{I}$ is isomorphic to the 
  quotient algebra of $S(S^{2}W)$ by the ideal generated by the image
  of $S^{2n+2}(W)$ under $P_{n+1}$.}}
\end{theorem}
{\bf Proof:} From Proposition \ref{invs} we see that $S^{2}W$ is
contained in $\mathcal{I}$. Consequently there is an algebra
homomorphism, $\phi$, from $S(S^{2}(W))$ to $\mathcal{I}$ which is
in fact $GL(W)$ equivariant. Indeed $\phi$ is
surjective which we know in any case
from Theorem 2 of \cite{howe}. To proceed we
use the following results,
\begin{enumerate}
\item $S(S^{2}W) = \oplus_{\mu} S_{\mu} W $ where $\mu$ is as in
  Proposition \ref{invs} except the second restriction is lifted,
  Proposition 1 in \cite{howe} page 562.\\
\item The ideal generated by $S^{2k}W$ is equal to $\oplus_{\mu}
  S_{\mu}(W)$ with $\mu$ as in Proposition \ref{invs} except that the
  second condition is replaced with the first row having $\geq 2k$
  elements, observing that the representation $S^{2k}W$ corresponds to
a Young diagram with only one row and it consists of $2k$ boxes. (A
special case of Theorem 3.1. of S. Abeasis \cite{abeasis})
\end{enumerate}
It follows that the ideal generated by $S^{2n+2}W$ is mapped to zero
under $\phi$ by using Proposition \ref{invs} (2) and the $GL(W)$
equivariance of $\phi$. Now using statement
1) above, Proposition \ref{invs} again and the surjectivity of $\phi$
we see that the ideal generated by $S^{2n+2}W$ is precisely the
kernel of $\phi$. This completes the proof.

\qed

\begin{remark}\label{P}{\rm{Let us write explicitly the map
      $P_{k}:S^{2k}W \rightarrow 
      S^{k}(S^{2}W)$. Fix a basis for $W$, $e_{1}, \dots , e_{g}$, let
      us also denote by $L_{ij}=e_{i}.e_{j}=L_{ji}$. Let
      $ e_{i_{1}}\dots
      e_{i_{2k}}$ with $i_{j} \in \{ 1, \dots , g\}$, be a monomial in
      $S^{2k}W$. Then 
      $P_{k}(e_{i_{1}}, \dots, e_{i_{2k}}) = $ the sum over all partitions in
      pairs $(i_{1},i_{2}),\dots (i_{2k-1}i_{2k})$ of the product
      $L_{i_{1},i_{2}}\dots L_{i_{2k-1}i_{2k}}$.}}
\end{remark}

\begin{definition}\label{irrep}{\rm{Given an integer $m$ and
      $\underline{\mu}=(\mu_{1}, \dots, \mu_{m} )$ with $\mu_{1} \geq
      \dots \geq \mu_{m} \geq 0$ we denote by $R_{\mu}(Sp(m))$ the
      irreducible representation of $Sp(m)$ with highest weight
      $\mu = \mu_{1}\lambda_{1} + \dots + \mu_{m}\lambda_{m}$, where
      the $\lambda_{i}$ are the fundamental weights. If it
      is clear which symplectic group one is referring to then we may
      write $R(\mu_{1}, \dots , \mu_{m})$ for the irreducible
      representation. 
}}
\end{definition}

\begin{theorem}\label{decomp}{\rm{We have the decomposition of $\Lambda(V \otimes
      W)$ under $Sp(V) \times Sp(g)$ as
$$
\Lambda(V \otimes W) = \oplus_{\mu} \,  R_{\tilde{\mu}}(Sp(V))\otimes
R_{\mu}(Sp(g))  
$$
where the allowed $\mu$ are $\underline{\mu} =
(n-a_{g}, \dots, n-a_{1})$, with $n \geq a_{1} \geq \dots \geq a_{g}
  \geq 0$ and $\tilde{\mu}$ is determined in terms of $\mu$ and given
  by
$$
\tilde{\mu} = \tilde{\omega}_{a_{1}} + \dots + \tilde{\omega}_{a_{g}}
$$
with $\tilde{\omega}_{i}$ is the $i$'th fundamental weight of
$Sp(V)$, alternatively one has
$$\tilde{\underline{\mu}} =( \sum_{i=1}^{g}\theta(a_{i}-1), \dots , \sum_{i=1}^{g}
 \theta(a_{i}-n)) 
$$
with
  $\theta(x) = 0$ if $ x< 0$ and $\theta(x)  = 1 $ if $x \geq
  0$. In particular, setting $\underline{\mu} = (n, \dots , n)$ so
  that $\underline{\tilde{\mu}}= 0$, we have
$$
\Lambda(V \otimes W)^{Sp(V)} = R_{(n,\dots , n)}(Sp(g)),
$$
that is the space of $Sp(V)$ invariants in $\Lambda(V \otimes W)$
carries the irreducible representation of $Sp(g)$ with highest weight
$(n , \dots , n)$.
}}
\end{theorem}
{\bf Proof:} This is implicit in the proof of Lemma 3.7 (especially pages
706-707) of \cite{BN}. However, note that the roles of the numbers
`$g$' and `$n$' are interchanged. As on p. 706 we find, for
$\underline{a}=(a_{1}, \dots , a_{g})$ with $n \geq a_{g} \geq \dots
\geq a_{1} \geq 0$, a vector $\phi_{\underline{a}} \in \Lambda (V
\otimes W)$ such that the 1-dimensional subspace spanned by
$\phi_{\underline{a}}$ is invariant under
$B_{1} \times B_{2}$ where $B_{1}$ and $B_{2}$ are suitable Borel
subgroups of $\mathrm{Sp}(V)$ and $\mathrm{Sp}(g)$ respectively. (It is
proven there that this subspace is invariant under $B_{1}\times B_{0}$
where $B_{0}$ is a Borel subgroup of $GL(W^{*})$ (indeed the upper
triangular matrices of $GL(W^{*})$) and also that
$\phi_{\underline{a}} \in \mathcal{H}$ in the notation of
\cite{BN}. Since $\mathcal{H}$ is the space annihilated by the
$L_{ij}$, defined in section \ref{explicit}, this in fact 
implies invariance under $B_{1} \times B_{2}$ as $B_{2}$ is the
subgroup generated by $B_{0}$ and the $L_{ij}$.)

The subrepresentation $\rho_{\underline{a}}$ spanned by the transforms
of $\phi_{\underline{a}} $ by $\mathrm{Sp}(V) \times \mathrm{Sp}(g)$
is irreducible with highest weight $\underline{\mu} = (n-a_{1}, \dots
, n-a_{g})$ for $\mathrm{Sp}(g)$ and $\tilde{\mu}$ for $\mathrm{Sp}(V)
$. As on page 707 of \cite{BN} we can show that the representations
$\rho_{\underline{a}}$ for $\underline{a}=(a_{1}, \dots , a_{g})$
satisfying  $n \geq a_{g} \geq \dots\geq a_{1} \geq 0$ exhaust all
irreducible representations contained in $\mathrm{Sp}(V) \times
\mathrm{Sp}(g)$. (Here we apply Theorem 8 of \cite{howe} to the
reductive dual pair $\mathrm{Sp}(V)$ and $\mathrm{Sp}(g)$ in the
orthogonal algebra.)

\qed

\begin{remark}{\rm{ By varying $W$ all irreducible representations of
      $\mathrm{Sp}(V)$ are realized on $\Lambda(V \otimes W)$.}}
\end{remark}

\begin{remark}{\rm{The case $g=1$ of the theorem is well known
      \cite{GW} (page 207 Corollary 4.5.9.).}}
\end{remark}

\begin{theorem}{\rm{
Let $F$ be a representation space of $\mathrm{Sp}(g)$ with
highest weight $(a_{1}, \dots , a_{g})$ such that $n \geq a_{1} \geq
\dots \geq a_{g} \geq 0$ then the ideal in
$S(S^{2}W)$ generated by the image of $P_{n+1}$
annihilates $F$.}}
\end{theorem}
\begin{remark}{\rm{Thus, as an
      ideal, the $P_{n+1}$ relations do not pick out a single
      irreducible representation, but rather, all those irreducible
      representations that occur in $\Lambda(V
      \otimes \mathbb{C}^{g})$.}}
\end{remark}

{\bf Proof:} Let $v\in F$ be a highest weight of some irreducible
representation $(a_{1}, \dots , a_{g})$ with $n \geq a_{1} \geq
\dots \geq a_{g} \geq 0$. 
We have at our disposal, by
Theorem  \ref{decomp}, all irreducible representations satisfying
$a_{1}\leq n$ (and 
ultimately all irreducible representations by taking symplectic vector
spaces $V$ of larger and larger dimension). As in section
\ref{explicit} the
$L_{ij}$ are realized as wedging with respect to the
symplectic form. In this case we may write
$$
P_{n+1}(e_{i_{1}},\dots , e_{i_{2n+2}}) = \mathcal{L}_{n+1}(I_{1},\dots, I_{2n+2})\,
v^{I_{1}}_{i_{1}} \wedge \dots \wedge v^{I_{2n+2}}_{i_{2n+2}}.
$$
As $P_{n+1}(e_{i_{1}},\dots, e_{i_{2n+2}})$ is totally symmetric in the
  labels $i_{1}, \dots , i_{2k}$, $\mathcal{L}_{n+1}(I_{1},\dots
  I_{2n+2})$ must be totally anti-symmetric in the labels $I_{1},\dots
  I_{2n+2}$, i.e. $\mathcal{L}_{n+1} \in \Lambda^{2n+2}V$ but as the
  dimension of $V$ is 2n there is no such tensor.

\qed

\section{Complex Manifolds}

Throughout this section $X$ is a compact closed complex manifold and
no extra structure is assumed. The dimension of $X$ is
$\dim_{\mathbb{C}}(X)=n$. Let $\Omega_{X}^{p}$
denote the sheaf of holomorphic $p$-forms on $X$. The objects of main
interest for us are the Dolbeault
cohomology groups of $X$ with values in
$$
\Omega_{X}^{\underline{p}} = \Omega_{X}^{p_{1}}\otimes \dots
\otimes \Omega_{X}^{p_{g}}
$$
where $\underline{p}$ is a $g$-tuple of natural numbers
$\underline{p}=(p_{1}, \dots, p_{g}) \in \mathbb{N}^{g}$, where
$\mathbb{N} = \{ 0, \, 1, \, 2, \, \dots \}$.

\begin{definition}{\rm{The pluri-Hodge groups of $X$ of
      length $g$ are
$$
\mathrm{H}^{q}(X,\Omega_{X}^{\underline{p}}) =
\mathrm{H}^{q}(X,\Omega_{X}^{p_{1}}\otimes \dots \otimes \Omega_{X}^{p_{g}} ).
$$}}
\end{definition}

\begin{remark}{\rm When $g-1$ of the $p_{i}$ are zero ($\Omega_{X}^{0} =
    \mathcal{O}_{X}$) the
    generalized Hodge groups are the Hodge cohomology groups
    $\mathrm{H}^{(p,q)}(X)$ .}
\end{remark}
\begin{definition}{\rm{The degree of $\alpha \in
    \mathrm{H}^{q}(X,\Omega_{X}^{\underline{p}})$ is $q + |p|$.}}
\end{definition}

\begin{definition}\label{dims}{\rm{The dimensions of the pluri-Hodge
    groups are
$$
h^{(\underline{p},q)}(X)= \dim{\mathrm{H}^{q}(X,
  \Omega_{X}^{\underline{p}})}.
$$}}
\end{definition}

\subsection{Pluri $\chi_{y}$ Genera determine Chern
  Numbers}\label{chernnumbers}

Now we motivate the introduction of pluri-Hodge groups by
showing that they determine the Chern numbers of $X$. To do this one
groups the dimensions $h^{(\underline{p},q)}(X)$ into a form
generalizing the usual $\chi_{y}$ genus. Then an application of
Riemann-Roch yields the result. 

\begin{definition}{\rm{The pluri-$\chi_{\underline{y}}$ genus
      (of length $g$) is,
$$
\chi_{-\underline{y}}(X) = \sum_{q, \; \underline{p}} (-1)^{q+
  |\underline{p}|} \, h^{(\underline{p},q)}(X)\, y_{1}^{p_{1}} \dots
y_{g}^{p_{g} }
$$
}}
\end{definition}
\begin{proposition}\label{expand}{{\rm
\be
\prod_{i=1}^{n}(a_{i} + te_{i}) = \sum_{p=0}^{n}t^{p}
\sum_{i_{1} < \dots < i_{p}} e_{i_{1}}\dots e_{i_{p}} \prod_{k \notin \{ i
\} } a_{k}.\nonumber
\ee
where $\{ i\}$ is the set $(i_{1},\dots , i_{p})$.
}}
\end{proposition}
No proof required.\\
\qed

\begin{proposition}\label{cclass}{{\rm Let $x_{r}$ denote the Chern
      roots of $X$ and $c_{i}$ the $i$-th Chern classes of the
      holomorphic tangent bundle of $X$. Chern classes are given by}
\be
c_{i} = \sum_{j_{1}< \dots < j_{n-i}} \frac{c_{n}}{x_{j_{1}}\dots
x_{j_{n-i}}} = \sum_{j_{1}< \dots < j_{n-i}} \prod_{k \notin \{j\}}
x_{k}. \nonumber
\ee
}
\end{proposition}
{\bf Proof:} We start with the definition in terms of Chern roots,
$$
\sum_{i=0}^{n}c_{i}\, t^{i} = \prod_{i=1}^{n} \left(1 +tx_{i}\right)
$$
so that
$$
\sum_{i=0}^{n}c_{i}\, t^{n-i} =  t^{n}\prod_{i=1}^{n}\left(1 + t^{-1} x_{i}
\right) = \prod_{i=1}^{n}\left(t + x_{i}\right)= c_{n} \,
\prod_{i=1}^{n}\left(tx_{i}^{-1} + 1\right) .
$$
We have, by proposition \ref{expand},
$$
\sum_{i=0}^{n}c_{i}\, t^{n-i} = c_{n}\,\sum_{p=0}^{n}t^{p}
\sum_{i_{1} < \dots < i_{p}} x_{i_{1}}^{-1}\dots x_{i_{p}}^{-1}.
$$
\qed

\begin{theorem}\label{cstring}{{\rm The coefficient of
      $(1-y_{1})^{n-q_{1}}\dots 
(1-y_{g})^{n-q_{g}}$ with $\sum_{i=1}^{g}q_{i}=\dim_{\mathbb{C}}(X)=n$ in the
      $\chi_{-\underline{y}}$
genus of length $g$ is $c_{q_{1}} \dots 
 c_{q_{g}}\, [X]$}.}
\end{theorem}
{\bf Proof:} The Riemann-Roch theorem, in the form presented in
\cite{NR}, tells us that
\be
\chi_{-\underline{ y} } = \mathrm{Todd} \, \left(\prod_{i=1}^{n} (a_{i} +
(1-y_{1})e_{i} \right) \dots \left(\prod_{i=1}^{n} (a_{i} +
(1-y_{g})e_{i} \right).[X],\nonumber
\ee
where, Todd is the Todd class of the tangent bundle of $X$ and
$$
e_{i} = e^{-x_{i}}, \;\;\; a_{i} = 1-e_{i}.
$$
Now from Proposition \ref{expand} we have that
\bea
\chi_{-\underline{ y} } & = & \mathrm{Todd}\, \sum_{\underline{p} } \left[
(1-y_{1})^{p_{1}} \sum_{i_{1} < \dots < i_{p_{1}}} e_{i_{1}}\dots
e_{i_{p_{1}}}  \prod_{k \notin \{ i
\} } a_{k} \dots \right. \nonumber \\
& & \;\;\; \left. \dots (1-y_{g})^{p_{g}}\sum_{j_{1} < \dots <
j_{p_{g}}} e_{j_{1}} \dots e_{i_{p_{g}}} \prod_{l \notin \{ j
\} } a_{l} \right] .[X]\nonumber
\eea
Set $p_{i} = n-q_{i}$, then the coefficient of $ (1-y_{1})^{n-q_{1}}
\dots (1-y_{g})^{n-q_{g}}$ is
\be
\mathrm{Todd}\,  \left[ \sum_{i_{1} < \dots < i_{n-q_{1}}}
e_{i_{1}}\dots e_{i_{n-q_{1}}} \prod_{k \notin \{ i 
\} } a_{k} \dots \sum_{j_{1} < \dots < j_{n-q_{g}}} e_{j_{1}}\dots
e_{i_{n-q_{g}}} \prod_{l \notin \{ j 
\} } a_{l} \right]. [X] .\label{coeff}
\ee

Note also that $a_{i} = x_{i} + \dots$ where
the ellipses represent forms of higher degree. Consequently,
\be
\prod_{k \notin \{ i
\} } a_{k} \dots \prod_{l \notin \{ j
\} } a_{l}, \label{1}
\ee
has form degree greater than or equal to $(n-p_{1}) + \dots (n-p_{g})=
\sum_{i=1}^{g}q_{i}$. Fix from now on $\sum q_{i} =n$. The product
(\ref{1}) is then a top degree form, and we may as well take
$a_{i}=x_{i}$, since all improvements yield an even higher degree
form. Since the products of the $a_{i}$ in (\ref{coeff}) already
constitute a form of highest possible degree we can set $\mathrm
{Todd} = 1$ and $e_{i}=1$. Hence (\ref{coeff}) simplifies to
\be
\left[ \sum_{i_{1} < \dots < i_{n-q_{1}}} \prod_{k \notin \{ i 
\} } x_{k} \dots \sum_{j_{1} < \dots < j_{n-q_{g}}} \prod_{l \notin \{ j 
\} } x_{l} \right] .[X]\nonumber
\ee
By Proposition \ref{cclass} this is
$$
c_{q_{1}}\, \dots \, c_{q_{g}}.[X]
$$
\qed

\begin{corollary}{{\rm Two complex manifolds of complex dimension $n$
are complex cobordant iff their generalized $\chi_{-\underline{y}}$
genera of length $n$ agree.}}
\end{corollary}
{\bf Proof:} Two complex manifolds are complex cobordant if all of
their Chern numbers agree. By Theorem \ref{cstring} all Chern
numbers are determined by the generalized $\chi_{-{\underline{y}}}$
genus of length $n$. Conversely the
Riemann-Roch theorem shows that the generalized
$\chi_{-\underline{y}}$ genera of any length are known once all Chern
numbers are.\\
\qed

\begin{remark}{\rm{The theorem is by no means sharp. For a complex
surface the $\chi_{y}$ genus determines and is determined by
the Chern numbers while the pluri-$\chi_{-\underline{y}}$ genus
of length 2 will do for up to
complex 5-folds. Also for a complex manifold with $c_{1}=0$, the 
pluri-$\chi_{y}$ genus of length $n/2$ suffices as the longest strings of
Chern numbers are either $c_{2}^{m}.[X]$ for $n=2m$ or $c_{2}^{m-1}c_{3}.[X]$
for $n=2m+1$.}} 
\end{remark}

\begin{remark}{\rm{While the $\chi_{-\underline{y}}$ genus tells us
      quite a bit about the pluri-Hodge groups it is not enough
      to establish that the dimensions are in fact deformation
      invariants or otherwise 
      of $X$. Clearly the combinations that appear in
      $\chi_{-\underline{y}}$ are invariant but we know that for
      non-K\"{a}hler manifolds even the usual Hodge numbers need not
      be invariants.}}
\end{remark}

\section{Holomorphic Symplectic Manifolds}\label{holsympman}

Let $X$ be a compact holomorphic symplectic
manifold. This means that $X$ is a compact complex manifold
with a nowhere degenerate holomorphic 2-form
$\eps$. Note that we do not assume that $\eps$ is closed. No extra
structure is assumed. In local coordinates write the components of
$\eps$ as $\eps_{IJ}$, since it is pointwise non-degenerate there
exists an inverse matrix at each point of $X$ which we denote by
$\eps^{IJ}$,
$$
\sum_{J} \eps^{IJ}\eps_{JK} = \delta^{I}_{K}.
$$

\subsection{Representations of the Symplectic Group on Sums of Pluri-Hodge
Groups}

Set
$$
\mathcal{H}^{q}(X) = \bigoplus_{\underline{p}} \mathrm{H}^{q}(X,
\Omega^{\underline{p} }_{X})= {\mathcal H}^{q}_{+}(X) \oplus {\mathcal
  H}^{q}_{-}(X)
$$
with
$$
{\mathcal H}^{q}_{+}(X) = \bigoplus_{|\underline{p}|= 0\; {\rm mod}\; 2}
\mathrm{H}^{q}(X,\Omega^{\underline{p}}_{X}), \;\;\;
{\mathcal H}^{q}_{-}(X) = \bigoplus_{|\underline{p}|= 1\; {\rm mod}\; 2}
\mathrm{H}^{q}(X,\Omega^{\underline{p}}_{X}).
$$
We shall now show that there is a natural action of $Sp(g)$ on the
holomorphic bundle $\Lambda( \Omega^{1}_{X} \otimes \mathbb{C}^{g})$
and hence on ${\mathcal H}^{q}_{\pm}(X)$ in case $X$ has a holomorphic
symplectic structure. 
\begin{theorem}\label{rep}{\rm{The holomorphic vector bundle $\Lambda(
      \Omega^{1}_{X} \otimes \mathbb{C}^{g})$ and the vector spaces
      ${\mathcal H}^{q}_{\pm}(X)$ 
      of pluri-Hodge groups of length $g$ form representations of
      ${\rm Sp}(g)$. 

}}
\end{theorem}
{\bf Proof:} The symplectic structure of $X$ gives a holomorphic
reduction of the cotangent bundle to the symplectic group
$\mathrm{Sp}(n)$ where $\dim_{\mathbb{C}} X = 2n$. Consider the action
of $\mathrm{Sp}(n)\times \mathrm{Sp}(g)$ on the vector space $\Lambda(
\mathbb{C}^{2n}\otimes \mathbb{C}^{g})$ discussed in section
\ref{sreps}. Since
the actions of $\mathrm{Sp}(n)$ and $\mathrm{Sp}(g)$ commute we obtain
an action of $\mathrm{Sp}(g)$, which is trivial on $X$, on the
associated vector bundle with typical fibre $\Lambda(
\mathbb{C}^{2n}\otimes \mathbb{C}^{g})$ which is the vector bundle $\Lambda(
      \Omega^{1}_{X} \otimes \mathbb{C}^{g})$ where $\mathbb{C}^{g}$ denotes
the trivial vector bundle. Hence we obtain an action of
$\mathrm{Sp}(g)$ on $\mathcal{H}^{q}(X, \Lambda(
      \Omega^{1}_{X} \otimes \mathbb{C}^{g})) =
      \mathcal{H}^{q}(X)$. Since this action preserves the degree we
      get an action on ${\mathcal H}^{q}_{\pm}(X)$.

\qed

\begin{remark}\label{remimp}{\rm{In other words: This is a direct
      application of the 
      construction of 
      section \ref{sreps} to 
obtain an $Sp(g)$ action on
$\Omega^{p_{1}}_{X} \otimes \dots \otimes \Omega^{p_{g}}_{X} $.
The only change to the
notation that is required is that $v_{i}^{I} \rightarrow dz_{i}^{I}$
where $dz_{i}^{I}$ are a local coordinate basis for $\Omega^{1}_{i}$
with $I = 1, \dots, n$ and $i=1, \dots, g$. We have therefore established
that the operators, $L_{ij}$, $\Lambda^{ij}$ and $h^{i}_{j}$ acting on
$\Omega^{p_{1}}_{X} \otimes \dots \otimes \Omega^{p_{g}}_{X}$ satisfy
the algebra of $sp(g)$ as given in Proposition \ref{spalg1} with
$\eps_{IJ}$ now being the components of the holomorphic symplectic
2-form. These operators commute with
$\overline{\partial}$, since they are made out the holomorphic 2-form, and
so the action extends to the pluri-Hodge groups. }}
\end{remark}

\begin{remark}{\rm{The $\mathrm{Sp}(g)$ action on $\alpha \in
\mathrm{H}^{q}(X,\Omega^{\underline{p}}_{X})$ is given by the formula
in Theorem \ref{rep1} with $v_{i}^{I} \rightarrow dz_{i}^{I}$ and
generalizes  the formulae with $n=1$ for the Hodge groups ($g=1$)
given in \cite{RW}.}}
\end{remark}

\begin{remark}{\rm{ Notice that the proof of the Lefschetz $sp(1)$ action
      with respect to the K\"{a}hler form is much more involved. This
      requires an application of Hodge theory since the corresponding
      Lie algebra generators do not all commute with $\overline{\partial}$.}}
\end{remark}

\subsection{The Geometric Meaning of the Genus}

In this section we assume that $X$ is compact. 
In \cite{T} it was shown that the $\chi_{y}$ genus is essentially the
same as the graded trace of an element of the $Sp(1)$ action on the
Hodge groups. Here this result is strengthened and the proof is
simplified. For ease of notation in the following theorem denote the condition
$|\underline{p}|=0\; \mathrm{mod}\, 2$ 
by $|\underline{p}|_{+}$ and the condition $|\underline{p}|=
1\; \mathrm{mod}\, 2$  by $|\underline{p}|_{-}$.
\begin{proposition}\label{prop}{\rm{Let X be a
      compact holomorphic  symplectic
  manifold of real 
dimension  $4n$. Let $U \in {\rm Sp}(g)$ and $y_{1}, \dots , y_{g},
  y_{1}^{-1}, \dots , y_{g}^{-1}$ be
eigenvalues of $U$, in the $2g$ dimensional defining representation,
then
$$
\Tr_{\mathcal{H}^{q}_{\pm}(X)}U =\frac{1}{y_{1}^{n}\dots y_{g}^{n}}
\sum_{|\underline{p}|_{\pm}} h^{\underline{p},q}(X)\,
y_{1}^{p_{1}}\dots y_{g}^{p_{g}}.
$$}}
\end{proposition}
{\bf Proof:} Since we have a representation the trace makes sense and
will give us a sum of characters of the group
element. Furthermore, since the space of diagonalizable
elements in $Sp(g)$ includes an open dense set, to establish the
formula we need only consider diagonalizable matrices. For
diagonalizable matrices, as the characters are class functions, one need
only consider their diagonalized form. Consequently, we may as well set 
$$U = \mathrm{diag} (y^{-1}_{1}, \dots,
y_{g}^{-1} ,y_{1},\dots,y_{g}).$$
By Theorem \ref{rep}, any form
$\alpha \in 
\mathrm{H}^{q}(X, \Omega_{X}^{\underline{p}})$ maps to
$$
y_{1}^{-n}\dots y^{-n}_{g}. y_{1}^{p_{1}}\dots y_{g}^{p_{g}}.\alpha
$$
so the action on the cohomology groups is just
multiplicative and we get such a factor for each element in
$\mathrm{H}^{q}(X, \Omega_{X}^{\underline{p}})$, that is we get such a
factor precisely $h^{(\underline{p},q)}(X)$ times.

\qed
\begin{remark}{\rm{The righthand side of the formula is invariant
      under sending any of the eigenvalues $y_{i}$ to their inverse
      $1/y_{i}$ since, by the triviality of the canonical bundle, we
      have $h^{(p_{1}, \dots, p_{i}, \dots,
      p_{g},q)}(X)=h^{(p_{1}, \dots, 2n-p_{i}, \dots, p_{g},q)}(X)  $.}}
\end{remark}

\begin{definition}{\rm{Let $U\in Sp(g)$, the graded trace or super
      trace on $\mathcal{H}(X)= \oplus_{q} \mathcal{H}^{q}(X)$ is
$$
\mathrm{STr}_{\mathcal{H}(X)} \, U = \sum_{q=0}^{2n}(-1)^{q}\left(
\Tr_{\mathcal{H}^{q}_{+}(X)}U - \Tr_{\mathcal{H}^{q}_{-}(X)}U\right).
$$}}
\end{definition}

\begin{theorem}\label{echi}{\rm{Let X be an irreducible holomorphic symplectic
  manifold of real 
dimension  $4n$. Let $U \in {\rm Sp}(g)$ and $y_{1}, \dots , y_{g}$ be
eigenvalues of $U$, in the $2g$ dimensional defining representation,
then
$$
{\rm STr}_{\mathcal{H}(X)} \, U = \frac{\chi_{-\underline{y}}}{y_{1}^{n} \dots
y_{g}^{n}}.
$$}}
\end{theorem}

{\bf Proof:} This is immediate given the previous proposition.

\qed

\subsection{Decomposition of $Sp(g)$ Representations on Pluri-Hodge Groups}

Fujiki \cite{F} gives a holomorphic version of the Lefschetz
decomposition theorem. This section is devoted to extending that
result to the pluri-Hodge groups.
\begin{definition}{\rm{Denote the multiplicity of the representation
      $R(a_{1}, \dots, a_{g})$ in $\mathcal{H}^{q}_{\pm}(X)$ by
      $m^{q}_{\pm}(a_{1}, \dots, a_{g})$.}}
\end{definition}

\begin{definition}{\rm{The Weyl group of $Sp(g)$ or the `octahedral'
      group of permutations and sign 
      changes of a $g$-tuple $(m_{1}, \dots , m_{g}) \in
      \mathbb{Z}^{g}$, as Weyl calls it, is denoted by $W_{g}$. For $w
      \in W_{g}$ set
      $[w]$ to equal the number of sign changes plus 1 if $w$ induces an
      odd permutation or simply equal to the number of sign changes if
      $w$ induces an even permutation.}}
\end{definition}

\begin{theorem}\label{multi}{\rm{Let $\underline{n} = (n, \dots, n)$,
      $\underline{\rho}= (g, \dots,1)$ and $\underline{a}= (a_{1}, \dots 
 , a_{g})$. We have, with the notation of Definition \ref{irrep} and
 Definition \ref{dims}, 
$$
m^{q}_{\pm}(a_{1}, \dots,
a_{g}) = \sum_{w \in W_{g}} (-1)^{[w]}\,
h^{(\underline{n}+w(\underline{a}+ \underline{\rho})-\underline{\rho}, q)}(X) .
$$}}
\end{theorem}

{\bf Proof:} The proof rests on an application of the decomposition
theorem of Weyl, 
(7.10.A) in \cite{W} (this is also known as outer multiplicity, see
Corollary 7.1.5 in \cite{GW}), combined with Proposition
\ref{prop} above. Since the proof is the same for either $\mathcal{H}_{+}^{q}$
and $\mathcal{H}_{-}^{q}$ we drop the distinction here. By Proposition
\ref{prop} and the remark after it we can write
$$
\Tr_{\mathcal{H}^{q}(X)}U =
\sum_{\underline{m}} h^{(\underline{n}+\underline{m},q)}(X)\,
y_{1}^{m_{1}}\dots y_{g}^{m_{g}}.
$$
According to the theorem of
Weyl, given a representation $V$ in which the character is a polynomial of
the form
$$
\chi_{V}(U)= \sum_{\underline{s}} k_{\underline{s}} \, y_{1}^{s_{1}} \dots
y_{g}^{s_{g}},
$$
the multiplicity of the representation $R(a_{1}, \dots, a_{g})$ is given
by
$$
m_{\underline{a}}=\sum_{w \in W_{g}} (-1)^{[w]}\,
k_{w(\underline{a}+ \underline{\rho})-\underline{\rho}} 
$$
on taking $V = \mathcal{H}^{q}(X)$ and comparing the two forms for the
character we have finished.\\
\qed

While the formula for the multiplicities is quite succinct it is
perhaps worth spelling it out for certain cases. The following
proposition is arrived at in the most pedestrian manner so the proof
is omitted.

\begin{proposition}\label{multi2}{\rm{The following multiplicities hold
\begin{enumerate}
\item For any $n$ and $g$ we have
$$
m^{q}(n, \dots, n, n-p)= h^{(p,q)}(X)-h^{(p-2,q)}(X),  \;\;\; 0\leq p
  \leq n .
$$
This includes the special case $g=1$
$$
m^{q}(n-p) =
  h^{(p,q)}(X)- h^{(p-2,q)}(X),
$$
which is one of the results of Fujiki.
\item We find for $Sp(2)$
\bea
m^{q}(n-p_{1},n-p_{2}) & = & h^{(p_{1},p_{2},q)}(X)-
h^{(p_{1},p_{2}-2,q)}(X) - h^{(p_{1}-4,p_{2},q)}(X)\nonumber\\
&+&
h^{(p_{1}-4,p_{2}-2,q)}(X)
 -  h^{(p_{2}+1,p_{1}-1,q)}(X) \nonumber\\
&&+ h^{(p_{2}+1,p_{1}-3,q)}(X) +
h^{(p_{2}-3,p_{1}-1,q)}(X) \nonumber \\
& & \;\;\; \; -h^{(p_{2}-3,p_{1}-3,q)}(X).\nonumber
\eea

\end{enumerate}
}}
\end{proposition}
\qed

\section{K3 Surfaces}

\subsection{The Pluri-Hodge Groups on a K3 Surface}

One expects, from physics, that the dimensions of the
pluri-Hodge groups for a holomorphic symplectic manifold $X$ are
deformation invariants of $X$. For $X=K3$ one can show this by
explicitly computing the dimensions. The following theorem is
provided by M.S. Narasimhan.

\begin{theorem}\label{k3constant}{\rm{Let $X$ be a K3 surface then the
      dimensions 
      $h^{(\underline{p},q)}(X) $ are independent of $X$. Set $m$ to
      equal the number of ones that appear in $\underline{p}$. The
      dimensions are given by, when $m$ is odd
$$
h^{(\underline{p},2)}(X) =h^{(\underline{p},0)}(X) = 0 \;\;\; {\rm and}\;\;\; 
h^{(\underline{p},1)}(X) = 2^{m+1}[6m - 1],
$$
while for $m$ even $m=2k$,
$$
h^{(\underline{p},2)}(X) =h^{(\underline{p},0)}(X)= \frac{(2k)!}{k!(k+1)!}
\;\;\; {\rm and}\;\;\;
h^{(\underline{p},1)}(X) = 2
\frac{(2k)!}{k!(k+1)!} + 2^{m+1}[6m-1] .
$$}}
\end{theorem}

{\bf Proof:} We will first show that the $h^{(\underline{p},0)}(X) $ =
constant. Choose a
K\"{a}hler metric on $X$ with Ricci curvature equal to zero. Recall that the
holonomy group is $SU(2)$. By Theorem 1.34 in
Kobayashi \cite{kob} every holomorphic
section of $\Omega^{1 \otimes m}_{X}$ is covariantly constant. 

Conversely every covariantly constant section of $(T^{*}X)^{\otimes m}$ defines
a holomorphic section of $(T^{*}X)^{\otimes m}$. In fact the connection
$D$ on $(T^{*}X)^{\otimes m}$ is the ``Chern connection'' on the
holomorphic bundle $(T^{*}X)^{\otimes m}$ so that if $D$ is the associated
covariant differential we have $D^{0,1} = \overline{\partial}$. Hence
if $D \alpha =0$, $\alpha$ a section of $(T^{*}X)^{\otimes m}$, then
$\overline{\partial} \alpha =0$.

Observe
that in our case the allowed values of $p_{i}$ are 0, 1, and
2. However, when $p_{i}=$ 0 or 2, $\Omega^{p_{i}}_{X}$ is the structure
sheaf and the theorem quoted above applies. Hence the
sections of the bundle are covariantly constant. Thus
$h^{(\underline{p},0)}(X)$ is equal to the dimension of the space of
invariants of $SU(2)$ in the m-fold tensor product of the defining
representation of $SU(2)$ where m is equal to the number of
$p_{i}=1$. This is clearly independent of $X$.\footnote{
More generally if $X$ is a compact Ricci flat K\"{a}hler manifold any
holomorphic section of $\Omega^{\underline{p}}_{X}$ is covariantly
constant and in principle can be 
calculated from the holonomy representations.}

As the canonical bundle is
trivial $\Omega^{p}_{X}$ is self dual and so by the duality theorem
we see that
$h^{(\underline{p},2)}(X)=h^{(\underline{p},0)}(X)$. Consequently
$h^{(\underline{p},1)}(X)$ is also constant. 

Now by calculating the dimension of the invariants and applying
the Riemann-Roch theorem we can calculate
$h^{(\underline{p},1)}(X)$. Riemann-Roch yields
$$
\int_{X} {\rm Todd}(X) \, {\rm ch}(\Omega^{1 \otimes m}_{X})= 2^{m+1}(1-6m),
$$
with $c_{2}(X) =24$.

On the other hand the space of invariants of $SU(2)$ on the m-fold
tensor product of the defining representation can be determined by
picking out the term in the character expansion of
\be
(e^{x}+e^{-x})^{m},
\ee
which corresponds to the trivial representation.
We want the dimensions of the $SU(2)$ invariants in $W^{\otimes m}$
where $W$ is the defining representation. Use the following fact:

(see e.g. \cite{S}) If ${\rm ch}(\lambda , V)$ is the character of a
representation of $SU(2)$ in a vector space $V$, then the dimension of
the invariants in $V$, (i.e. the multiplicity of the trivial
representation) is given by the co-efficient of $e^{\lambda}$ in
$(e^{\lambda}-e^{-\lambda}){\rm ch}(V)$.

We have ${\rm ch}(W) = e^{\lambda}+ e^{-\lambda}$ and ${\rm
ch}(W^{\otimes m}) = (e^{\lambda}+ e^{-\lambda})^{m}$, hence
\bea
(e^{\lambda}-e^{-\lambda}){\rm ch}(W^{\otimes m}) &=&
(e^{\lambda}-e^{-\lambda}) (e^{\lambda}+ e^{-\lambda})^{m} \nonumber
\\
& = & \sum_{r=0}
\left( \begin{array}{c}
m \\
r
\end{array} \right) e^{\lambda(m+1 -2r)} - \sum_{s=0}
\left( \begin{array}{c}
m \\
s
\end{array} \right) e^{\lambda(m -1 -2s)} \nonumber.
\eea
Now we can get $e^{\lambda}$ either from $m =2r$ or $m = 2s + 2$ in
either case $m$ is necessarily even. Thus if $m$ is odd the trivial
representation does not occur in $W^{\otimes m}$. On the other hand
when $m=2k$ the coefficient of $e^{\lambda}$ is
$$
\left( \begin{array}{c}
2k \\
k
\end{array} \right) - \left( \begin{array}{c}
2k \\
k-1
\end{array} \right) = \frac{(2k)!}{k!(k+1)!} .
$$
\qed

\section{Graph Representations of $Sp(g)$}\label{graphs}

We saw in previous sections that the $P_{n+1}$ relations determine
certain representations of $Sp(g)$. It is quite remarkable that the
same relations appear in the construction of the LMO invariant. They
are imposed as a set of relations so that the second Kirby move be
respected. 

In this section a certain amount of pictorial gymnastics is
unavoidable. For a more complete account of chord diagrams and the
$IHX$, $AS$ (and orientation), $STU$, $P_{n}$, $O_{n}$ and branching
relations $B$, I refer the reader to \cite{LMO, MO}.

\begin{definition}{\rm {A chain graph, $\Gamma_{g}$, with $g$ copies
    of $S^{1}$, is depicted in Figure \ref{chaingraph}. 
\begin{figure}[h]
\centerline{\epsffile{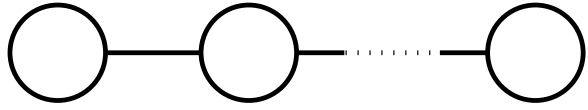}}
\caption{The chain graph $\Gamma_{g}$.}\label{chaingraph}
\end{figure}

A chord diagram $D$ with support $\Gamma_{g}$ is $\Gamma_{g}$
    together with an oriented uni-trivalent graph with the univalent
    vertices lying in $\Gamma_{g}$. Denote by
    $\mathcal{C}(\Gamma_{g})$ the
    $\mathbb{C}$-linear space of chord diagrams with support
    $\Gamma_{g}$ including trivial circles
    of dashed lines.}} 
\end{definition}

\begin{definition}{\rm{Denote the number of 
trivalent vertices by
      $D^{i}$ (i-for internal) and let $D^{i}_{>n}$ be the relation
      such that any chord diagram with internal degree more than $n$
      is zero. Let $D_{T}$, the degree, be half the sum of trivalent
      and univalent vertices. 
}}
\end{definition}

\begin{definition}{\rm{We set 
$\mathcal{A}(\Gamma_{g})
    = \mathcal{C}(\Gamma_{g})/ IHX, AS, STU, B$ and furthermore let 
$\mathcal{A}^{(n)}(\Gamma_{g}) =$\\
 $ 
      \mathcal{A}(\Gamma_{g})/O_{n},\, P_{n+1}, \,
      D^{i}_{>2n}$. }}
\end{definition}
\begin{remark}{\rm{The $STU$ relations do not preserve the internal
      degree $D^{i}$ even though they do preserve the total degree
      $D_{T}$ which
      counts trivalent as well as univalent vertices. We will need a
      slightly different presentation of chord diagrams which does not
      make use of the $STU$ relations. Such a space of uni-trivalent graphs
      was introduced by Sawon \cite{justinhk}. 
}}
\end{remark}

\begin{definition}{\rm{A marked uni-trivalent graph is an oriented
      uni-trivalent graph whose univalent vertices carry labels $1$ to
      $g$ and repetitions are allowed. Denote by $\mathcal{B}_{g}$ the
      space of linear combinations of marked uni-trivalent graphs
      modulo $AS$ and $IHX$ and including trivial circles
    of dashed lines. We also denote the completion of
      $\mathcal{B}_{g}$ with 
      respect to the total degree by the same symbol.
}}
\end{definition}
\begin{remark}{\rm{Sawon did not include dashed circles in his
      definition, we will need this generalization below.
}}
\end{remark}
\begin{proposition}{\rm{\cite{justinhk} There is an isomorphism
      $\tau:\mathcal{B}_{g} \rightarrow \mathcal{A}(\Gamma_{g})$.
      
}}
\end{proposition}
\begin{remark}{\rm{The isomorphism requires that $\mathcal{B}_{g}$
      include the `empty' uni-trivalent graph, $e_{g}$, which corresponds to
      $\Gamma_{g} \in \mathcal{A}(\Gamma_{g})$. Furthermore, it is not
      difficult to show that the isomorphism
      between $\mathcal{B}_{g}$ and $\mathcal{A}(\Gamma_{g})$ given in
      \cite{justinhk} preserves the $P_{n}$ and $O_{n}$ equivalence relations.
}}
\end{remark}
We now see that there is a natural $Sp(g)$ action on these spaces.
\begin{theorem}\label{Brep}{\rm{The space of graphs $\mathcal{B}_{g}/O_{n}$ is
      a (reducible) representation space of $Sp(g)$.}}
\end{theorem}
{\bf Proof:} The manner in which the Lie algebra of $Sp(g)$ is to be
realised on the space of graphs $\mathcal{B}_{g}/  O_{n}$ is as follows:

$L_{ij}$ introduces a dashed line connecting two univalent vertices one with
the label $i$ and the other with the label $j$.
$h^{i}_{\; j}$ relabels an existing univalent vertex with the label
$i$ with a $j$  
and one repeats this summing over all univalent vertices with the
label $i$. The
action of $\Lambda^{ij}$ 
is best described as follows. Consider that $i\neq j$ first. If there
are no univalent vertices marked either $i$ or $j$ the action of
$\Lambda^{ij}$ is to give zero. So suppose that there are such marked univalent
vertices. For fixed such marked univalent vertices
$\Lambda^{ij}$ acts to eliminate both the univalent vertices and it joins the
dashed lines that ended on those vertices and multiplies by
$-1/4$. The action of $\Lambda^{ij}$ on the uni-trivalent graph is to
sum over the operation just described on all possible pairs of
univalent vertices with the markings $i$ and $j$. $\Lambda^{ii}$ is
defined in the same way except that the multiplicative factor is
$-1/2$.
\begin{figure}[h]
\centerline{\epsffile{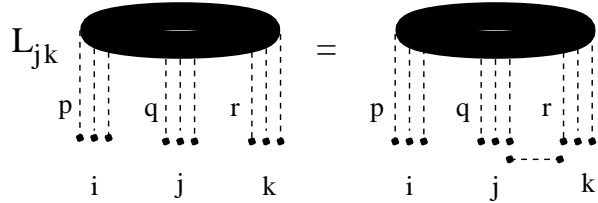}}
\caption{$L_{ij}$ acting on a uni-trivalent graph.}\label{Algebra1fig}
\end{figure}

\begin{figure}[h]
\centerline{\epsffile{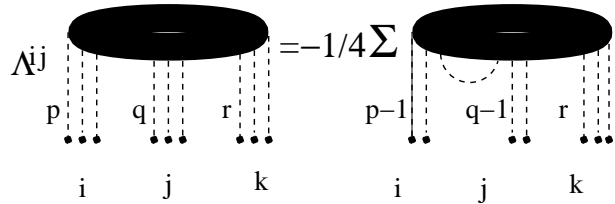}}
\caption{$\Lambda^{ik}$ acting on a uni-trivalent graph with the sum being over
  all pairs of univalent vertices marked $i$ and $j$.}\label{Algebra2fig} 
\end{figure}

\begin{figure}[h]
\centerline{\epsffile{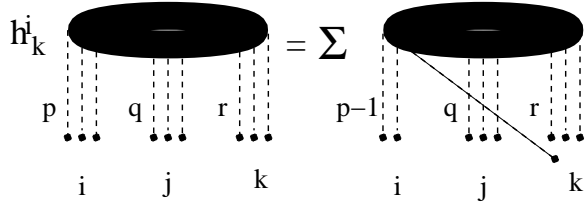}}
\caption{The action of $h^{i}_{k}$ on a uni-trivalent graph, the sum is over
the $p$ univalent vertices marked $i$.}\label{Algebra3fig}
\end{figure}

One can now show graphically that these maps satisfy the algebra given in
Proposition \ref{spalg1}. For example in Figures \ref{Algebra1fig},
\ref{Algebra2fig} and 
\ref{Algebra3fig} the filled in ellipse designates the rest of a
uni-trivalent graph, the labels $p$, $q$
and $r$ indicate the number of vertical dashed lines connected to
univalent vertices marked with $i$, $j$ and $k$ respectively and if a
dashed line is not strictly 
vertical it is not counted. One
checks that, on these figures, $[\Lambda^{ij}, L_{jk}] = -1/4
h^{i}_{k}$, $[h^{i}_{\; j}, h^{j}_{\; k}] = - h^{i}_{\; k}$ and so
on. In this manner, running through the various possibilities, one 
establishes the claim that the space of uni-trivalent graphs
$\mathcal{B}_{g} /  O_{n}$ is a representation space of $Sp(g)$.

\qed
\begin{definition}{\rm{Let the vector space
    $\mathcal{B}^{n}_{g \, m} $ denote the part of
      $\mathcal{B}_{g}/P_{n+1},O_{n}$ with internal degree $m$, and set
      $\mathcal{B}^{(n)}_{g} =
      \mathcal{B}_{g}/P_{n+1},O_{n}, D^{i}_{>2n}$ so that
      $\mathcal{B}^{(n)}_{g} =\sum_{m=0}^{2n}\oplus
      \mathcal{B}^{n}_{g\, m}  $.
      }}
\end{definition}

\begin{figure}[h]
\centerline{\epsffile{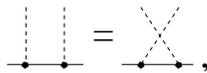}}
\caption{The TU Relation.}\label{TUfig}
\end{figure}

\begin{lemma}\label{mlemma}{\rm{Let $D$ be a uni-trivalent graph in
      $\mathcal{B}_{g}$. Then, modulo the equivalence relation
      $P_{n+1}$, 
\begin{enumerate}
\item The uni-trivalent
  graph $D$ is equivalent to a linear sum of uni-trivalent graphs
  which have at most $2n$ univalent vertices labeled $i$ for all $i
  \in \{ 1, \dots , g\}$ and
\item The uni-trivalent graph $D$ becomes equivalent to a linear sum
  of uni-trivalent graphs each of which has either $n$ dashed chords
  joining vertices with the same labels or at most $2n-1$ univalent
  vertices with a given marking $i$.
\end{enumerate}
}}
\end{lemma}
{\bf Proof:} This is a transcription to $\mathcal{B}_{g}$ of Lemma 3.1
in \cite{LMO}. The proof there can be adapted to the current situation
as follows. All vertices with the same marking can be thought of as
being connected to the component $C$ in that paper. Since ordering of
univalent vertices with the same marking is
immaterial in $\mathcal{B}_{g}$ the univalent vertices on $C$ obey a $TU$
relation, see Figure \ref{TUfig}, rather than the $STU$
relation. Consequently, the
only difference with \cite{LMO} is that now we do not
generate `lower order terms' and the proof goes through.

\begin{figure}[h]
\centerline{\epsffile{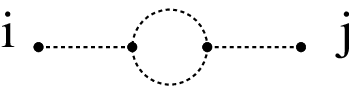}}
\caption{Two marked univalent vertices joined to two trivalent
  veritices}\label{sym0} 
\end{figure}

\begin{figure}[h]
\centerline{\epsffile{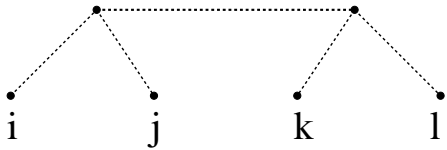}}
\caption{Four marked univalent vertices joined to two trivalent
  vertices}\label{g4} 
\end{figure}

\begin{figure}[h]
\centerline{\epsffile{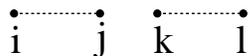}}
\caption{A pair of chords joining marked univalent vertices.}\label{g24}
\end{figure}

\begin{figure}[h]
\centerline{\epsffile{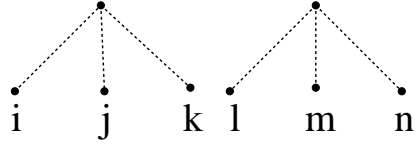}}
\caption{Six marked univalent vertices joined to two trivalent
  vertices}\label{sym} 
\end{figure}

\begin{figure}[h]
\centerline{\epsffile{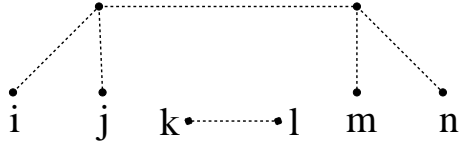}}
\caption{Equivalence with previous graphs.}\label{g6mod}
\end{figure}

\qed
\begin{proposition}{\rm{$\mathcal{B}_{g\, 2}^{1} \cong \theta \sqcup
      \mathcal{B}_{g\, 0}^{1} $ where $\sqcup$ is disjoint union and
      $\theta$ is the graph 
\begin{figure}[h]
\centerline{\epsffile{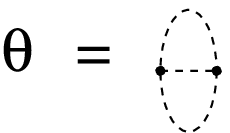}}
\end{figure} 
}}
\end{proposition}
{\bf Proof:} There are two trivalent vertices. If they are not
connected to any univalent vertex then they form a $\theta$. Now
suppose that the trivalent vertices do not form a $\theta$ 
diagram. There are a number of different situations to
consider. By Lemma
\ref{mlemma} we only need consider that if a univalent
vertex marked $i$ is connected to a trivalent vertex then, it is the unique
such marked univalent vertex. Consider the
situation where the two trivalent vertices are joined by two dashed lines,
($g \geq 2$)
as in Figure \ref{sym0}.
An application of $P_{2}$ and IHX shows that such a graph is
equivalent to -1/2
$\theta\, \sqcup $ a dashed chord joining the $i$ and $j$ univalent
vertices. Two trivalent vertices can be 
connected by just one dashed line if $g \geq 4$, as shown in Figure \ref{g4}.
An application of $P_{2}$ and IHX turns such uni-trivalent graphs into sums
of graphs of the form of $\theta$ union chord diagrams as in Figure
\ref{g24}. 
Lastly if $g\geq 6$ we can have uni-trivalent graphs where the two trivalent
vertices are not connected at all by a dashed line as in
Figure \ref{sym}.
By repeated use of $P_{2}$ such a uni-trivalent graph can be expressed as
a sum of terms of the form shown in Figure \ref{g6mod}. However, these
graphs have already been dealt with successfully and so we are
done. 

\qed

\begin{remark}{\rm{From the work of J. Murakami \cite{Mu} we see that
$\mathcal{A}^{(1)}(\Gamma_{1})$ and
$\mathcal{A}^{(1)}(\Gamma_{2})$ are two and five-dimensional spaces as 
$\mathcal{A}^{(1)}(\emptyset)$-modules. We note that $\dim R(1)
=2$ and that $\dim R(1,1)=5$ where $R$ is a representation of
$Sp(g)$. We are now in a position to show that
this numerical agreement is no accident.}}
\end{remark}

\begin{theorem}\label{chordspg}{\rm{The space of uni-trivalent graphs
      $\mathcal{B}^{(1)}_{g}$ is a representation space for
      $Sp(g)$. In fact, as $Sp(g)$ representation spaces,
      $\mathcal{B}^{n}_{g\, 0} 
      = R(n, \dots , n)$ and $\mathcal{B}^{n}_{g\, 1}= R(n,
      \dots , n,
      0,0,0)$.

}}
\end{theorem}
{\bf Proof:} 
By the previous theorem $\mathcal{B}_{g}/ O_{n}$ is a representation
space for $Sp(g)$, we will now see that the $P_{n+1}$ equivalence
relations pick out certain irreducible representations. 

For $\mathcal{B}^{n}_{g\, 0}$ any uni-trivalent graph can
be obtained from the empty graph, $e_{g}$, by successive
application of $L_{ij}$ for various $i$, $j$. Consequently, all the
uni-trivalent graphs belong to one irreducible representation of
$Sp(g)$ since, for a semi-simple Lie algebra, a finite dimensional reducible
representation decomposes into irreducible representations. The
representation is finite dimensional since the $P_{n+1}$ relations
imply that no more than $2n$ 
univalent vertices have the same marking, $D^{i}=0$ says that there
are no trivalent 
vertices and $O_{n}$ ensures that dashed circles are set to be
multiplication by $-2n$. 

$e_{g}$ is a highest weight state with
weight $(n, \dots , n)$. To see this we apply Theorem 4.7.3 in
\cite{V}. To adapt to the notation of that reference let $Y_{g}
= L_{gg}$, $X_{g}= \Lambda^{gg}$ and $Y_{i} = h^{i+1}_{i}$ 
and $X_{i}= h^{i}_{\; i+1}$ for $i =1, \dots, g-1$ and set $H_{i}= n-
h^{i}_{i}$. Clearly $X_{i} e_{g}=0$, $H_{i}e_{g}= n e_{g}$
and $Y_{i}^{n+1} e_{g}=0$, proving that $e_{g}$ is a highest
weight vector with the advertized weight.

In the case of $\mathcal{B}^{n}_{g \, 1}$, by Lemma \ref{mlemma},
every uni-trivalent graph can be 
obtained as follows. Begin by connecting by dashed lines the unique
trivalent vertex to three distinctly marked univalent vertices (if any
two of the univalent vertices have the same marking $AS$ assures us that this
vanishes). Then add an
even number of univalent vertices (of any marking) to the graph and join these
in pairs by dashed lines. 

Let $D^{1}_{g}(i,j,k)$ be the uni-trivalent graph consisting of a
trivalent vertex connected by dashed lines to three univalent vertices
with markings $i$, $j$ and $k$. Acting with $L_{mn}$ one generates
part of the span of 
$\mathcal{B}^{n}_{g\, 1}$ just described
and by summing over all possible $i$, $j$ and $k$
one obtains the complete span of $\mathcal{B}^{n}_{g\, 1}$. However, as $
D^{1}_{g}(i,j,k)= h^{g-2}_{\; i}h^{g-1}_{\; j}h^{g}_{\;
  k}D^{1}_{g}(g-2,g-1,g)$ one can obtain any uni-trivalent graph in
$\mathcal{B}^{n}_{g\, 1}$ starting
from $D^{1}_{g}(g-2,g-1,g)$. This uni-trivalent graph is a highest weight
with $\lambda(H_{g})= \lambda(H_{g-1})=\lambda(H_{g-2})=0$ and
$\lambda(H_{i}) = n$. 

\qed

\section{Rozansky-Witten Topological Field Theory}

I now want to make contact between the pluri-Hodge groups, the
graph representations of the previous section and so of the
surprising omnipresence of the skew representation theory of
$Sp(g)$. The pluri-Hodge groups and uni-trivalent graphs are brought
together in the context of the theory of Rozansky and Witten \cite{RW}.
\begin{remark}{\rm{On any complex compact closed manifold, $Y$,
the Atiyah class \cite{At} $\alpha \in
\mathrm{H}^{1}(Y,\Omega^{1}_{Y} \otimes \Omega^{1}_{Y} \otimes T_{Y}) $
measures the obstruction to having a global holomorphic section
of the holomorphic tangent bundle $T_{Y}$ of $Y$. When $Y=X$ is
holomorphic-symplectic there is an identification $T_{X} \cong
\Omega^{1}_{X}$ (through the holomorphic 2-form) and we also denote by
$\alpha$ its image in
$\mathrm{H}^{1}(X, T_{X} \otimes T_{X} \otimes T_{X}) $. Indeed $
\alpha \in \mathrm{H}^{1}(X, \mathrm{Sym}^{3}T_{X}) $. 
}}
\end{remark}
\begin{proposition}{\rm{\cite{justinhk} Let $X$ be a compact closed
       holomorphic symplectic manifold with $\dim_{\mathbb{C}}X=2n$. There is a
      homomorphism of vector spaces 
$$
W_{X}: \mathcal{B}_{g}/O_{n} \rightarrow \mathcal{H}_{g}(X).
$$
If a marked uni-trivalent graph $D$ has $q$ trivalent vertices and
$p_{i}$ univalent vertices marked by $i$ then
$$
W_{X}(D) \in \mathrm{H}^{q}(X, \Omega_{X}^{p_{1}} \otimes \dots
\otimes \Omega_{X}^{p_{g}}).
$$
}}
\end{proposition}
\begin{definition}{\rm{$W_{X}$ is defined as follows: Suppose there
      are $q$ trivalent vertices in the graph. To each such vertex we
      associate the Atiyah class $\alpha$, and cup product is
      understood between the Atiyah 
      classes. Each dashed leg of the vertex is understood to
      correspond to one of the fibre factors $T_{X}$ in $T_{X}\otimes
      T_{X} \otimes T_{X}$. If there is a dashed line between two
      trivalent vertices, then
      one contracts the fibre factors of the contracted legs of the
      two different trivalent vertices as
      follows $T_{X} \otimes T_{X} 
      \cong T_{X}\otimes \Omega^{1}_{X} \rightarrow
      \mathcal{O}_{X}$. If a trivalent vertex is connected to a
      univalent vertex marked $i$ then the $T_{X}$ factor of the leg
      attached to the trivalent vertex is replaced with the $i$'th
      $\Omega^{1}_{X}$. If two univalent vertices, marked $i$ and $j$
      respectively, are connected by a dashed line then this is mapped to
      $\eps_{ij}$ and the cup product is understood
      throughout. $W_{X}(e_{g})= 1 \in \mathrm{H}^{0}(X,
      \mathcal{O}_{X})$ and the circle of a dashed line maps to $-2n$.
}}
\end{definition}
\begin{remark}\label{nextr}{\rm{That $W_{X}$ is compatible with the
      $AS$ and $IHX$ 
      relations when there are no univalent vertices was shown in
      \cite{RW}, but perhaps a more accessible reference is
      \cite{HS}. That the image of $W_{X}$ is in $\mathcal{H}_{g}(X)$
      as well as why $W_{X}$ is compatible with the $AS$ and $IHX$
      relations when there are marked univalent vertices
      is explained in some detail in \cite{justinhk}. Essentially, in
      local coordinates we have
      $$ W_{X}(D) = \sum F_{i_{1} \dots j_{p_{g}}} dz_{1}^{i_{1}}
      \wedge \dots \wedge 
      dz_{1}^{i_{p_{1}}} \otimes \dots \otimes dz_{g}^{j_{1}} \wedge
      \dots \wedge 
      dz_{g}^{j_{p_{g}}}$$
where the coefficients $F_{i_{1} \dots j_{p_{g}}}$ (which are
themselves $(q,0)$-forms) are
totally anti-symmetric in all labels not just within each marking.
}}
\end{remark}
\begin{theorem}\label{commspg}{\rm{The weight system
      $W_{X}:\mathcal{B}_{g}/O_{n}  \rightarrow
\mathcal{H}_{g}(X) $ commutes with the respective $Sp(g)$ actions
on $\mathcal{B}_{g}/O_{n}$ and $\mathcal{H}_{g}(X)  $.
}}
\end{theorem}
{\bf Proof:} Let $D$ be some uni-trivalent graph, then by the
definition of $W_{X}$ we have that
$W_{X}(L_{ij} D) = \eps_{ij} \wedge W_{X}(D)$. However, wedging with
$\eps_{ij}$ is the 
same as acting with the $Sp(g)$ 
generators $L_{ij}$ of section \ref{holsympman}, see Remark
\ref{remimp}. It is not difficult to see, given Remark \ref{nextr},
that the rest of the 
generators map correctly. Compatibility with the $O_{n}$ relation is a
property of 
the holomorphic 2-form, namely that $\sum_{I,
  \,J}\eps^{IJ}\,\eps_{JI}= -2n$. \\
\qed 
\begin{corollary}\label{cork3}{\rm{Let $X$ be a K3 surface, then
      $W_{X}:\mathcal{B}^{1}_{g \, 0} \rightarrow
      \mathcal{H}_{+}^{0}(X) $ and $W_{X}:\theta \sqcup
      \mathcal{B}^{1}_{g \, 0} \rightarrow \mathcal{H}^{2}_{+}(X)$ and,
     for $g\geq 3$, $W_{X}:\mathcal{B}^{1}_{g \, 1}\rightarrow
      \mathcal{H}^{1}_{-}(X)$.
}}
\end{corollary}
{\bf Proof:} This follows from the theorem and Theorem
\ref{chordspg}. Note that $W_{X}(\theta) \in \mathrm{H}^{2}(X,
\mathcal{O}_{X})$ a non-zero multiple of the
anti-holomorphic 2-form.\\
\qed
\begin{remark}{\rm{The multiplicity of the trivial representation in
$\mathcal{H}^{q}(K3)$ with $g=3$ can be deduced from Theorem \ref{multi} to
be $h^{(1,1,1,q)}(K3)-2h^{(1,q)}(K3)$ and by Theorem
\ref{k3constant} this is $232$ for $q=1$ and zero otherwise. Of this
large number of possibilities it is $\mathbb{C}.\alpha$ that is the
trivial $Sp(3)$ representation of the corollary and more generally
$\alpha \in \mathrm{H}^{1}(X, \Omega_{X}^{\underline{p}})$ with
$p_{g}=p_{g-1}=p_{g-2}=1$ and all other $p_{i}=0$ is the highest
weight state of $R(1, \dots , 1 ,0,0,0)$ for $g \geq 3$. 
}}
\end{remark}

\begin{remark}{\rm{Sawon \cite{justinhk} noted that $W_{X}$ is not
      surjective (as is evident from Corollary \ref{cork3}). Indeed
      we will now see that for any holomorphic 
      symplectic manifold $X$
      it misses at least `1/2' of the pluri-Hodge groups.
}}
\end{remark}
\begin{proposition}\label{half}{\rm{The pluri-Hodge groups $\mathrm{H}^{q}(X,
      \Omega_{X}^{\underline{p}}) $ are not in the image of $W_{X}$
      for $q+ |\underline{p}| \not\in 2\mathbb{N}$ or for
      $q=|\underline{p}|=1$.
}}
\end{proposition}
{\bf Proof:} If there are $q$ trivalent vertices in a graph $D$ these
have $3q$ legs. Each univalent vertex can be thought of as a `leg' so
there are $|\underline{p}|$ of these. All legs must be joined in pairs
so that $3q + |\underline{p}|$ must be even to have such a graph in
$\mathcal{B}_{g}$. When $q=1$ the $AS$ relation ensures that $D$ is
equivalent to zero if there is only one univalent vertex.\\
\qed

\begin{remark}{\rm{The chord diagrams, $\mathcal{A}(\Gamma_{g})$, of
      the previous section arise in the topological quantum field
      theory (TQFT) of Murakami and Ohtsuki. Sawon introduced $W_{X}$
      as a `hyper-K\"{a}hler weight system' for the MO
      invariants. This is the map $W_{X} \circ \tau^{-1}:
      \mathcal{A}(\Gamma_{g}) \rightarrow \mathcal{H}_{g}(X)$. The
      natural hope being that this would be the correct TQFT approach
      to the Rozansky-Witten theory \cite{RW}. One can show that this
      expectation is borne out for $X$ a $K3$ surface. Let $M$ be a
rational homology 3-sphere and $K$ a null-homologous knot in
$M$. Denote the complement of a tubular neighbourhood of $K$ in $M$ by
$M \backslash K$. Rozansky and Witten determine (using path integrals)
the associated vector in $\mathcal{H}_{1}(K3)$; this is (5.68) in
\cite{RW} and denoted by $| M\backslash K \rangle$ there. The Murakami
invariant \cite{Muinv}, $\Lambda_{1}(M,K)$, is the first MO
invariant in a special normalization. For $K$ with zero
framing we have
$W_{K3}\circ \tau^{-1}(\Lambda_{1}(M,K)) = | M\backslash K \rangle$
\cite{Tunpub} (though the formulae for both $| M\backslash K \rangle$ and
$\Lambda_{1}(M,K)$, given in \cite{RW} and \cite{Muinv} respectively,
require some minor corrections).
}}
\end{remark}

\begin{remark}{\rm{Physics implies the much more interesting result
      that the pluri-Hodge groups provide a representation of the
      mapping class group $\mathrm{MCG}_{g}$ of a Riemann surface
      $\Sigma_{g}$ of genus $g$. J. Murakami \cite{Mu} gives 
      representations of $\mathrm{MCG}_{g}$ on
      $\mathcal{A}^{(1)}(\Gamma_{g})$ for  
      $g=1$ and $g=2$ with a non-trivial Torelli subgroup action. The
      map $W_{K3}\circ \tau^{-1}$ induces representations of
      $\mathrm{MCG}_{g}$ on 
      $\mathcal{H}^{0}_{+}(K3) \oplus \mathcal{H}^{2}_{+}(K3) $ (which
      mixes the two spaces), with a non-trivial Torelli action
      \cite{justinmapclass, NT}.
}}
\end{remark}

\rnc{\Large}{\normalsize}

\end{document}